\documentclass[12pt,3p,times]{elsarticle}
\usepackage{amsmath,amssymb,amsfonts}
\usepackage{algorithmic}
\usepackage{graphicx}
\usepackage{textcomp}
\usepackage{xcolor}
\usepackage{caption}
\usepackage{subcaption}
\usepackage{hyperref}
\usepackage{natbib}
\usepackage{float}
\usepackage{array}
\usepackage{makecell}
\usepackage{booktabs}
\usepackage{lipsum}
\usepackage{siunitx}

\newtheorem{thm}{Théorème}[section]

\begin{document}

\begin{frontmatter}

\title{A Novel Hyperchaotic System Derived from Asymmetric Bidirectional Coupling of Two Lorenz Oscillators: Synchronization Analysis and Cryptographic Applications}

\author[m1]{Mimoun Hamdi}
\ead{my.hamdi@gmail.com}

\author[m1]{Mohamed Karim Hamdani\corref{cor1}}
\ead{hamdanikarim42@gmail.com}
\cortext[cor1]{CONTACT\; M.K. Hamdani: hamdanikarim42@gmail.com}

\author[ipeik1,ipeik2]{Kamel Abboudi}
\ead{kamalo1982@yahoo.fr}

\begin{center}
%\address[m0]{Military Academy of Tunis, Military Academy Foundouk Jedid, 8012 Nabeul, Tunisia}
\address[m1]{Science and Technology for Defense Laboratory, Center for Military Research l'Aouina, Tunisia}
\address[ipeik1]{Preparatory Institute for Engineering Studies, University of Kairouan, Kairouan, Tunisia}
\address[ipeik2]{Laboratory of Modeling, Mechanics and Production Engineering, National Engineering School of Sfax (ENIS), Sfax, Tunisia}
\end{center}

\begin{abstract}
This paper investigates a six-dimensional coupled Lorenz system obtained through bidirectional asymmetric diffusive coupling of two identical Lorenz oscillators. The proposed model is analyzed from both dynamical and synchronization perspectives. Dissipativity conditions, equilibrium points, local stability properties, and Lyapunov spectra are systematically examined. The results reveal wide regions of chaotic and hyperchaotic behavior, with up to three positive Lyapunov exponents depending on the coupling parameters. These findings provide a practical framework for selecting coupling configurations that guarantee hyperchaotic dynamics while preserving the simplicity of the classical Lorenz structure.
\\
A Lyapunov-based synchronization analysis is then developed. By exploiting effective attractor bounds obtained from the actual system trajectories rather than conservative worst-case estimates, a less restrictive sufficient synchronization condition is derived. Under the classical Lorenz parameters, synchronization is guaranteed for $\alpha_1+\alpha_2>15.1$.
\\
The synchronization performance is further evaluated through a parametric analysis of the synchronization time. Numerical results show that increasing the coupling strength significantly accelerates convergence toward the synchronization manifold. For a broad range of admissible coupling values, synchronization times below $0.5$ second are achieved, whereas previously reported Lorenz-based synchronization schemes typically require approximately $5$ seconds.
\\
Finally, a secure communication framework based on synchronized hyperchaotic generators is presented. Image encryption and decryption experiments confirm the effectiveness of the proposed approach. The combination of controllable hyperchaotic dynamics, explicit synchronization conditions, and fast synchronization performance makes the proposed coupled Lorenz system attractive for secure communication and chaos-based cryptographic applications.
\end{abstract}


\end{frontmatter}

\section{Introduction}
Chaotic systems have evolved from theoretical nonlinear dynamical models into powerful tools for a wide range of engineering applications, including secure communications, cryptography, signal processing, random number generation, and information security \cite{pecora1990}. Among the numerous chaotic models reported in the literature, the Lorenz system remains one of the most studied and influential examples due to its simple structure, rich nonlinear behavior, and extreme sensitivity to initial conditions \cite{lorenz1963}. Its well-known butterfly-shaped attractor has made it a benchmark model for the investigation of chaos generation, synchronization, and chaos-based communication schemes.
\\
As illustrated in Figure \ref{fig:chaotic_trajectories}, the classical Lorenz system generates complex chaotic trajectories and can be synchronized under appropriate coupling or control strategies. However, the increasing security requirements of modern communication systems have motivated the search for dynamical systems exhibiting higher complexity than that provided by conventional three-dimensional chaotic oscillators. In particular, single chaotic systems possess only one positive Lyapunov exponent, which may limit their unpredictability, effective key space, and resistance against reconstruction or parameter estimation attacks in cryptographic applications \cite{li2005, kanso2012}.
\\
To overcome these limitations, considerable attention has been devoted to hyperchaotic systems. Unlike conventional chaotic systems, hyperchaotic systems possess at least two positive Lyapunov exponents, leading to multiple independent expanding directions in phase space. As a result, they exhibit higher dynamical complexity, stronger sensitivity to initial conditions and parameters, improved mixing properties, and larger key spaces, making them particularly attractive for secure communication and encryption applications \cite{Biswas2026, gao2008, grassi2002, Ogabi2025, Zhou2026}.
\\
Several approaches have been proposed to generate hyperchaotic behavior, including dimensional augmentation, nonlinear feedback mechanisms, and the coupling of multiple chaotic oscillators \cite{chen2004, wang2008}. Among these approaches, coupled chaotic systems offer an effective framework for increasing dynamical complexity while preserving relatively simple mathematical structures. Numerous studies have investigated synchronization and control strategies for coupled Lorenz systems, including adaptive control, nonlinear feedback control, sliding-mode approaches, and observer-based synchronization techniques \cite{bragard2007, grassi2008}. Despite these efforts, the derivation of explicit synchronization conditions remains challenging due to the strong nonlinearities involved, and many existing results rely on conservative assumptions or complex controller designs.

\begin{figure}[H]
\centering
\includegraphics[width=15cm,height=7cm]{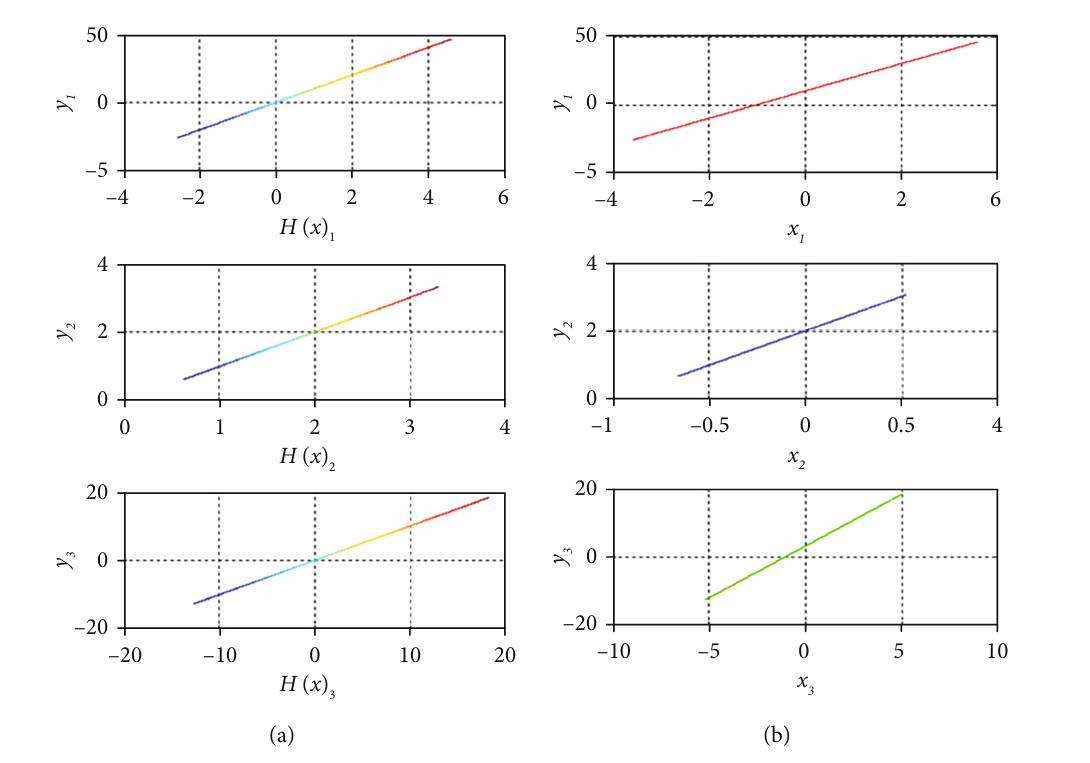}
\caption{Synchronization relations between state variables in the unidirectionally coupled Chen system. (a) Complete synchronization on the manifold \( H(x) \); (b) generalized synchronization relative to \( y = H(x) \).}
\label{e2022}
\end{figure}

Motivated by these observations, this work investigates a six-dimensional coupled Lorenz system obtained through bidirectional asymmetric diffusive coupling. Rather than introducing a fundamentally new chaotic model, the objective is to provide a comprehensive dynamical and synchronization analysis of this coupled configuration and to evaluate its suitability for secure communication applications. The proposed coupling structure employs two independent coupling coefficients, allowing the investigation of a broad range of interaction scenarios and synchronization regimes.
\\
The scientific contribution of this work lies in the combination of theoretical analysis, numerical investigation, and application-oriented validation. First, the dissipativity properties of the proposed system are established, and the admissible parameter region ensuring phase-space contraction is determined. Second, the equilibrium points and their local stability properties are analyzed through Jacobian eigenvalue analysis. Third, the Lyapunov spectrum is investigated over a wide range of coupling parameters, allowing the identification of chaotic and hyperchaotic operating regions. Fourth, sufficient synchronization conditions are derived using a Lyapunov-based framework combined with Sylvester's criterion. In contrast to classical worst-case approaches based on overly conservative global bounds, the proposed analysis exploits effective attractor bounds obtained numerically from the system trajectories, leading to less conservative synchronization conditions. Furthermore, the synchronization performance is evaluated through a synchronization-time analysis and compared with results reported in the literature.  We refer the reader to \cite{Akter2023,Chauhan2023,Ogabi2025} for an overview of references on this subject. In 2015, Khan and Singh~\cite{KS2015} studied the synchronization of the Lorenz chaotic system using a nonlinear control strategy. They considered the following  {master system}:
\begin{align}\label{eq:drive}
\begin{cases}
\dot{x}_1(t) = \sigma \big(x_2(t) - x_1(t)\big), \\
\dot{x}_2(t) = r x_1(t) - x_2(t) - x_1(t) x_3(t), \\
\dot{x}_3(t) = -b x_3(t) + x_1(t) x_2(t),
\end{cases}
\end{align}
where \( \sigma \), \( r \), and \( b \) are positive system parameters.  
\\
The  {slave system} (or response system) is defined as:
\begin{align}\label{eq:response}
\begin{cases}
\dot{y}_1(t) = \sigma \big(y_2(t) - y_1(t)\big) + u_1(t), \\
\dot{y}_2(t) = r y_1(t) - y_2(t) - y_1(t) y_3(t) + u_2(t), \\
\dot{y}_3(t) = -b y_3(t) + y_1(t) y_2(t) + u_3(t),
\end{cases}
\end{align}
where \( u_i(t),\, i = 1, 2, 3 \), are the control inputs. The errors between the master and slave states are defined as
\[
e_i(t) = y_i(t) - x_i(t), \quad i = 1,2,3.
\]
The control laws are chosen as:
\[
u_1(t) = k_{11} e_1(t), \quad
u_2(t) = k_{22} e_2(t), \quad
u_3(t) = k_{33} e_3(t),
\]
with the assumption that the master system states are bounded:
\[
|x_1(t)| \leq M_1, \quad |x_2(t)| \leq M_2, \quad |x_3(t)| \leq M_3.
\]
The error dynamics are then given by:
\begin{equation}\label{eq:error}
\begin{cases}
\dot{e}_1(t) = \sigma e_2(t) - (\sigma + k_{11}) e_1(t), \\
\dot{e}_2(t) = (r - M_3) e_1(t) - (1 + k_{22}) e_2(t) - M_1 e_3(t) + e_1(t) e_3(t), \\
\dot{e}_3(t) = M_2 e_1(t) + M_1 e_2(t) - (b + k_{33}) e_3(t) - e_1(t) e_2(t).
\end{cases}
\end{equation}
Khan and Singh~\cite{KS2015} established the following theorem:
\begin{thm}
The master system~\eqref{eq:drive} and the slave system~\eqref{eq:response} are synchronized if the control gain matrix \( B = \mathrm{diag}(k_{11}, k_{22}, k_{33}) \) is chosen such that the error system~\eqref{eq:error} is asymptotically stable. This is ensured if there exists a positive definite Lyapunov function
\[
V(e,t) = e^T P e
\]
such that its time derivative satisfies
\[
\dot{V}(e,t) = e^T Q e < 0 \quad \text{for all } e \neq 0.
\]
\end{thm}
For the specific control functions
\[
u_1(t) = 0, \quad u_2(t) = k_{22} e_2(t), \quad u_3(t) = k_{33} e_3(t),
\]
and the conditions:
\begin{enumerate}
\item \( \sigma > 0 \),
\item \( k_{22} < \frac{1}{4}(\sigma + r - M_3)^2 - 1 \),
\item \( k_{33} < \frac{M_2^2 (1 + k_{22})}{4\sigma (1 + k_{22}) - (\sigma + r - M_3)^2} - b\),
\end{enumerate}
they showed that \( \dot{V}(e) < 0 \) for all \( e \neq 0 \), which implies the errors converge to zero, hence synchronization. Numerical simulations validated these results using:
\[
M_1 = 12, \quad M_2 = 15, \quad M_3 = 48, \quad r = 28, \quad b = \frac{8}{3},
\]
and gains:
\[
\sigma = 20, \quad k_{22} = -1, \quad k_{33} = -\frac{8}{3}.
\]
Initial conditions for the chaotic master and slave systems were:
\[
x(0) = [10, -10, 20], \quad y(0) = [7, -10, 10],
\]
and for synchronized systems:
\[
x(0) = [50, -10, 20], \quad y(0) = [20, 10, 20].
\]
As shown in Figure \ref{fig:chaotic_trajectories}, the Lorenz system exhibits complex chaotic trajectories and can be effectively synchronized between two identical systems under appropriate control.
\begin{figure}[H]
\centering
\includegraphics[width=8cm,height=3cm]{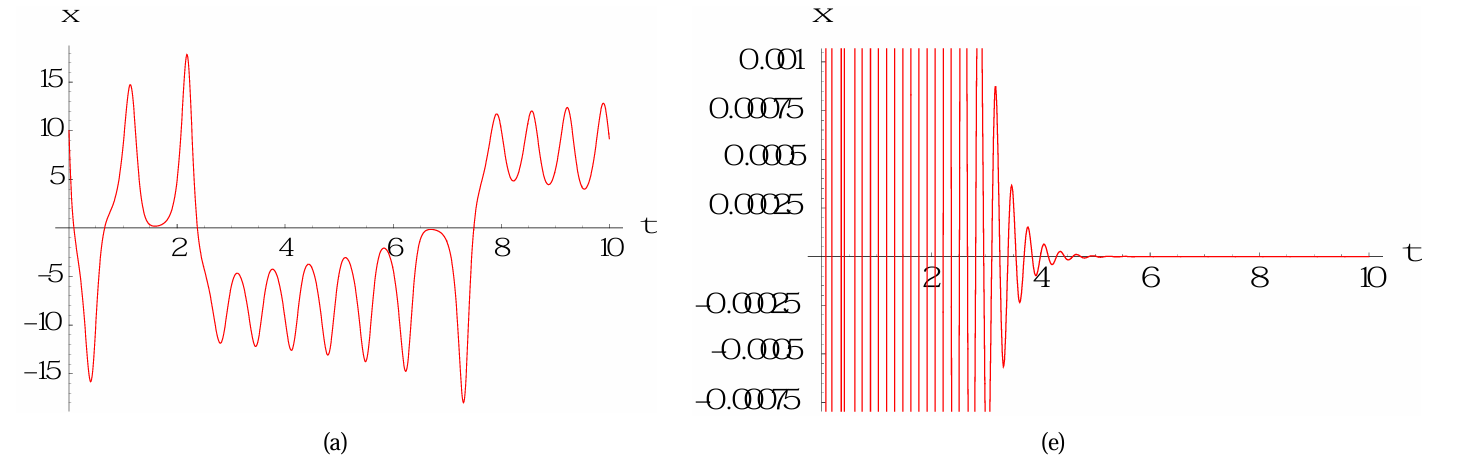}
\includegraphics[width=8cm,height=3cm]{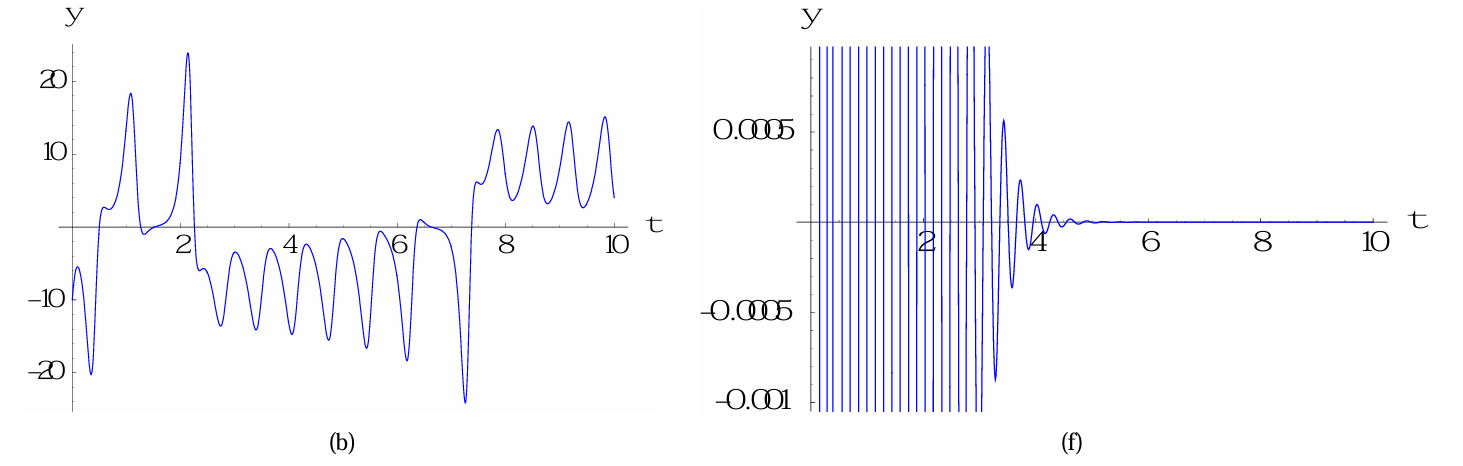}
\end{figure}
\begin{figure}[H]
\includegraphics[width=8cm,height=3cm]{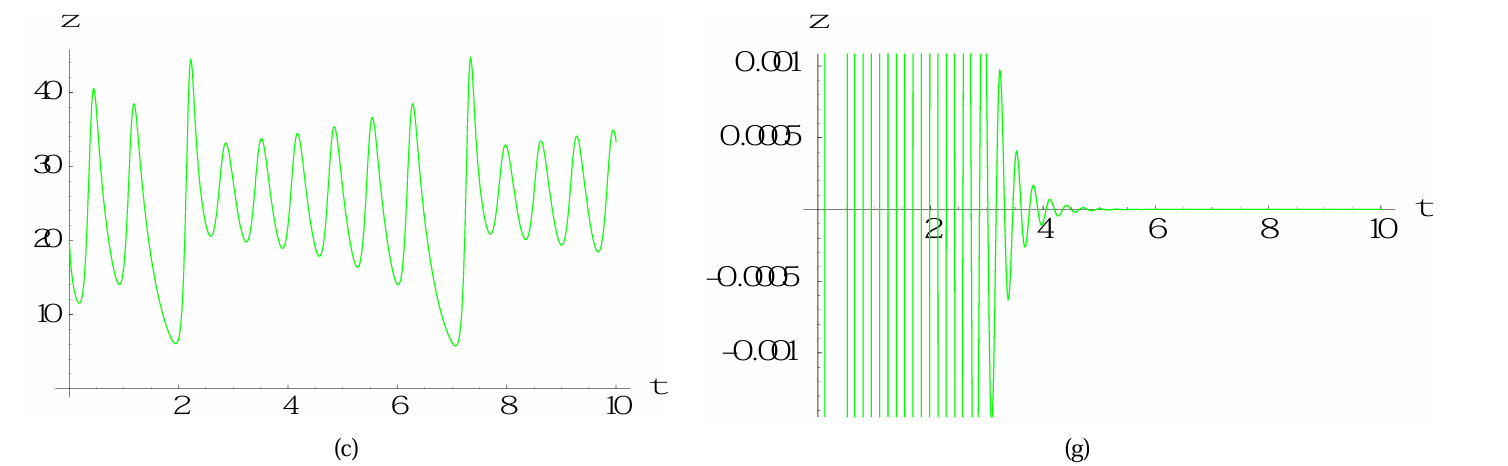}
\includegraphics[width=8cm,height=3cm]{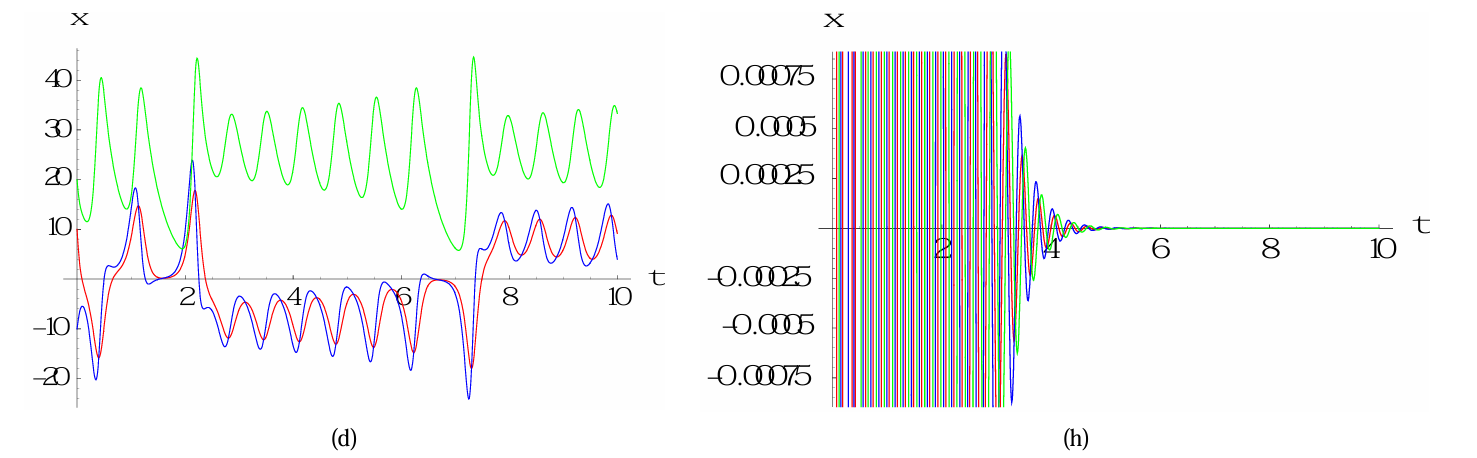}
\caption{The first column (a)--(d) shows chaotic trajectories of the Lorenz system. The second column illustrates effective synchronization between two identical systems under control.}
\label{fig:chaotic_trajectories}
\end{figure}
\noindent
Recently, in 2022, Yu, Zuo, Zhu, and Zhang \cite{YZZZ2022} studied two unidirectionally coupled dynamical systems:
\begin{equation}\label{tag1}
\dot{x} = F(x),
\end{equation}
\begin{equation*}\label{tag2}
\dot{y} = G(x, y),
\end{equation*}
where \( x \in \mathbb{R}^n \), \( y \in \mathbb{R}^m \), \( F(x) \) is a differentiable vector field, and \( G(x,y) \) describes the response system dynamics, depending on both \( x \) and \( y \).
\\
For the chaotic system~\eqref{tag1}, the response system is constructed using a linear error-feedback approach:
\begin{equation}\label{tag3}
\dot{y} = D H(x) \cdot F(x) + K \left( y - H(x) \right),
\end{equation}
where \( D H(x) \) is the Jacobian of \( H(x) \) and \( K = \mathrm{diag}(k_1, k_2, \dots, k_m) \) is a diagonal coupling matrix.
\\
The authors proposed a general criterion for \textit{unidirectional generalized chaotic synchronization} relative to a target manifold \( y = H(x) \).
\begin{thm}[Yu, Zuo, Zhu, and Zhang~\cite{YZZZ2022}]
If the coupling matrix \( K \) is chosen such that
\[
\lambda_i < 0, \quad i = 1, 2, \dots, m,
\]
where \( \lambda_i \) are eigenvalues of \( K^T P + P K \) for a positive definite symmetric matrix \( P \), then the unidirectionally coupled systems~\eqref{tag1} and~\eqref{tag3} achieve generalized synchronization along the manifold \( y = H(x) \).
\end{thm}
An application was performed for the Chen chaotic system~\cite{CU1999}:
\begin{equation}\label{eq:chen}
\begin{cases}
\dot{x}_1 = a(x_2 - x_1), \\
\dot{x}_2 = (c - a)x_1 - x_1 x_3 + c x_2, \\
\dot{x}_3 = x_1 x_2 - b x_3,
\end{cases}
\end{equation}
with \( a = 35, b = 3, c = 28 \). The target function is \( H(x) = [x_1, e^{x_2}, e^{-x_3}]^T \), initial conditions are
\[
x(0) = (-10,0,27)^T, \quad y(0) = H(x(0)) = (-9,2,84)^T,
\]
and the diagonal coupling matrix is
\[
K = \mathrm{diag}(-0.5, -0.5, -0.5).
\]
The synchronization behavior in unidirectionally coupled systems is illustrated in Figure \ref{e2022}, showing the relationship between state variables in the unidirectionally coupled Chen system.

\noindent
Yu and Zhang~\cite{YZ2004} studied bidirectionally coupled Lorenz systems:
\begin{equation}\label{eqYZ:coupled_lorenz}
\begin{cases}
\dot{x}_1 = a(x_2 - x_1) + d_{11}(y_1 - x_1), \\
\dot{x}_2 = c x_1 - x_1 x_3 - x_2 + d_{12}(y_2 - x_2), \\
\dot{x}_3 = x_1 x_2 - b x_3 + d_{13}(y_3 - x_3), \\
\dot{y}_1 = a(y_2 - y_1) + d_{21}(x_1 - y_1), \\
\dot{y}_2 = c y_1 - y_1 y_3 - y_2 + d_{22}(x_2 - y_2), \\
\dot{y}_3 = y_1 y_2 - b y_3 + d_{23}(x_3 - y_3),
\end{cases}
\end{equation}
with \( a = 10, b = 8/3, c = 28 \). Without coupling (\( d_{ij} = 0 \)), each system evolves independently and exhibits chaotic behavior (Figure \ref{fig:chaotic_trajectories}).

\begin{figure}[H]
\centering
\includegraphics[width=9cm,height=3cm]{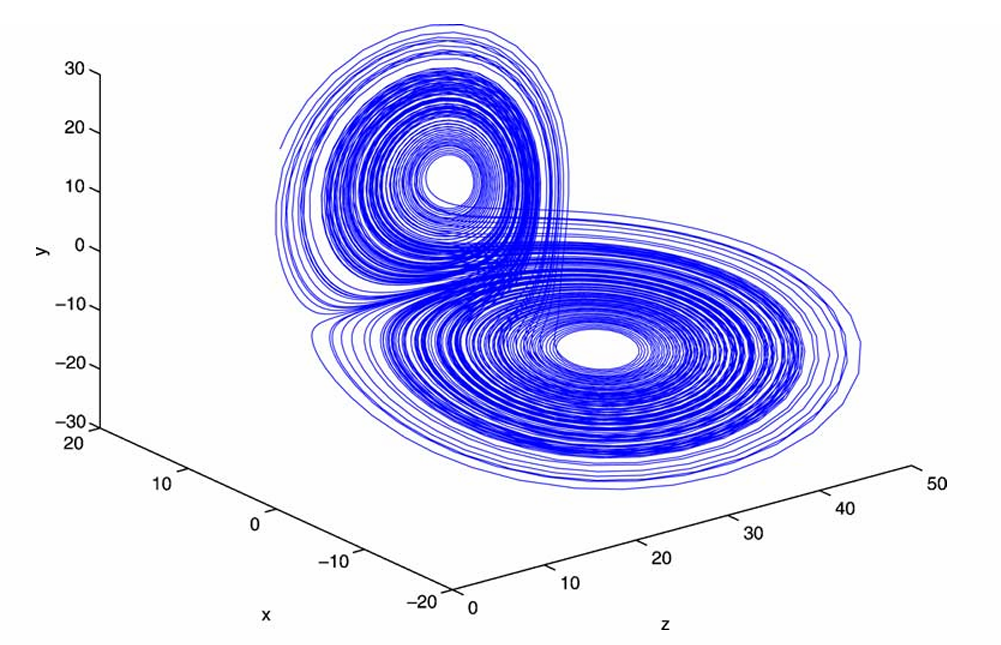}
\caption{Chaotic attractor of the Lorenz system.}
\label{fig:chaotic_trajectories}
\end{figure}

\begin{figure}[H]
\centering
\begin{tabular}{ccc}
\textbf{(a)} & \textbf{(b)} & \textbf{(c)} \\
\includegraphics[width=5.1cm,height=3cm]{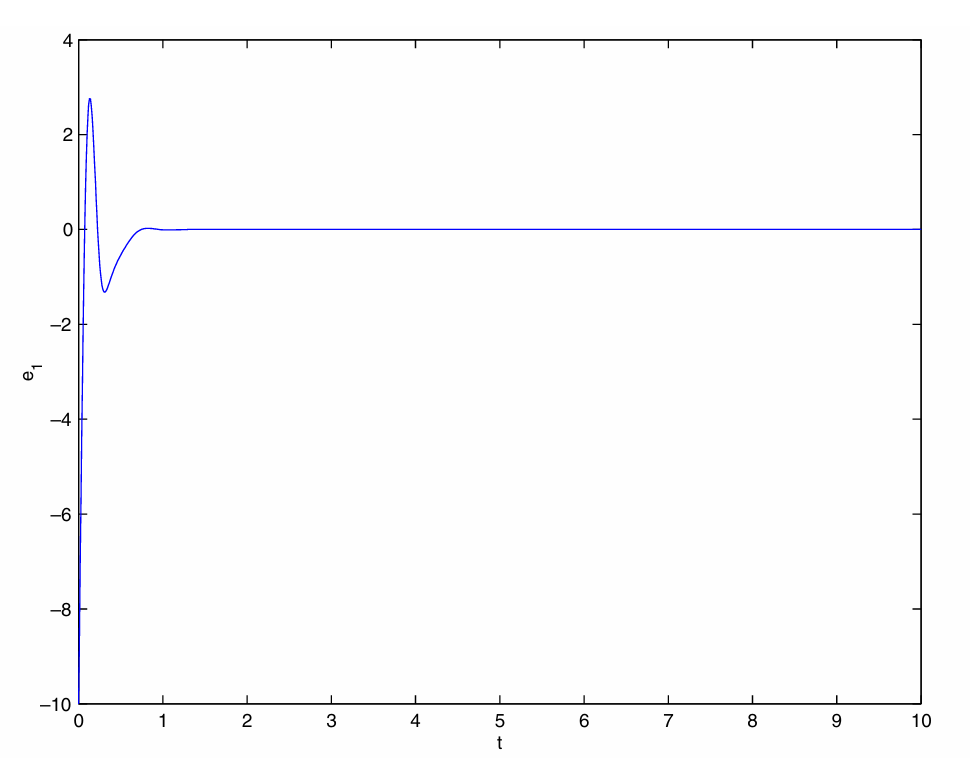} &
\includegraphics[width=5.1cm,height=3cm]{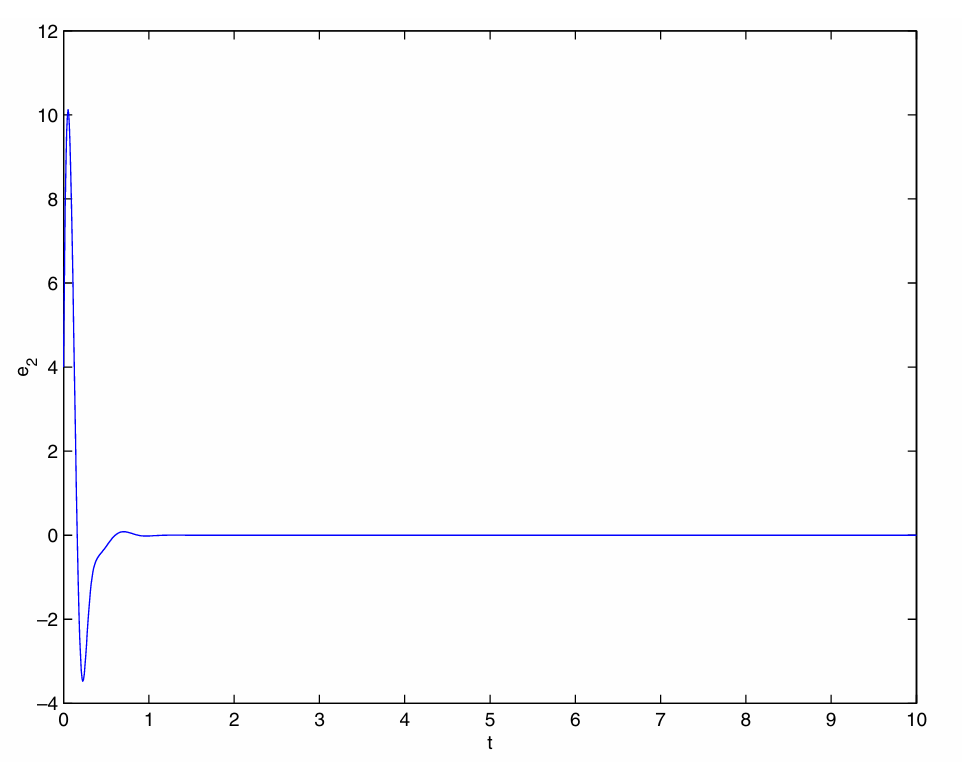} &
\includegraphics[width=5.1cm,height=3cm]{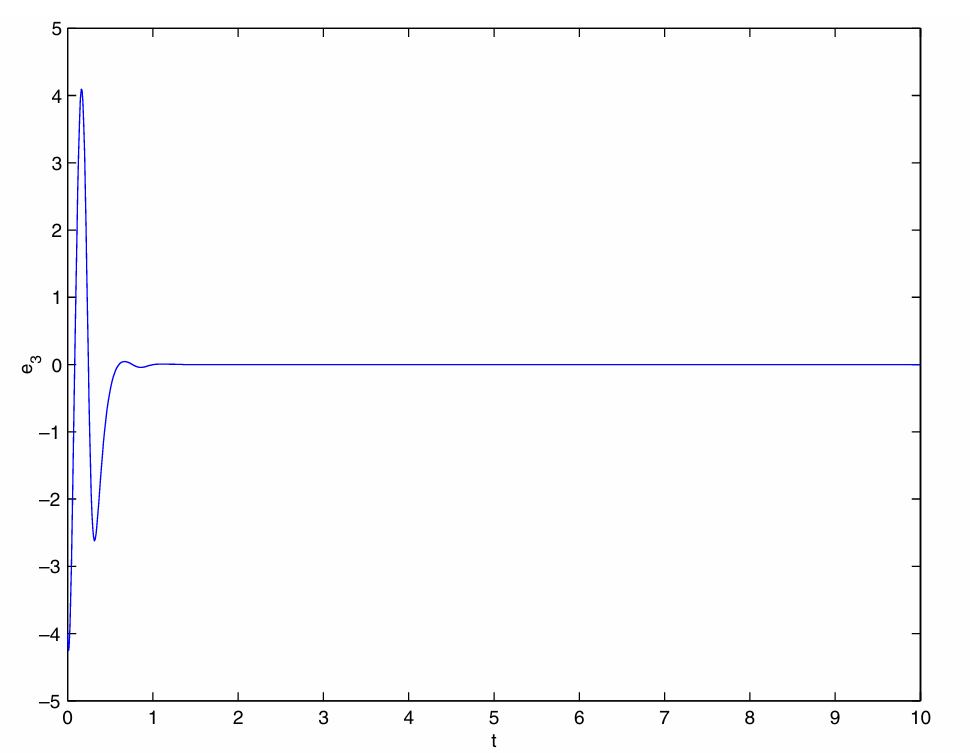} \\
\end{tabular}
\caption{Synchronization errors of the bidirectionally coupled Lorenz systems: (a) evolution of \(e_1=x_1-y_1\), (b) evolution of \(e_2=x_2-y_2\), and (c) evolution of \(e_3=x_3-y_3\).}
\label{fig:errors}
\end{figure}

Their criterion ensures global synchronization if the coupling parameters satisfy some conditions and a positive definite matrix \( P \) exists. Numerical simulations illustrate the synchronization errors \( e_1, e_2, e_3 \) as shown in Figure \ref{fig:errors}.
Motivated by these observations, this work investigates a novel six-dimensional dynamical system generated through asymmetric bidirectional diffusive coupling of two identical Lorenz oscillators. Unlike conventional symmetric coupling structures, the proposed model introduces two independent coupling coefficients, $\alpha_1$ and $\alpha_2$, allowing non-reciprocal interactions between the subsystems. Such a framework provides greater flexibility for modeling realistic communication and networked systems while simultaneously enriching the resulting dynamics. The proposed system comprises two identical Lorenz oscillators coupled through asymmetric diffusive terms, mathematically represented as:
\begin{equation}
\begin{cases}
\dot{x}_1 = \sigma(y_1 - x_1) + \alpha_1(x_2 - x_1) \\
\dot{y}_1 = x_1(\rho - z_1) - y_1 + \alpha_1(y_2 - y_1) \\
\dot{z}_1 = x_1y_1 - \beta z_1 + \alpha_1(z_2 - z_1) \\
\dot{x}_2 = \sigma(y_2 - x_2) + \alpha_2(x_1 - x_2) \\
\dot{y}_2 = x_2(\rho - z_2) - y_2 + \alpha_2(y_1 - y_2) \\
\dot{z}_2 = x_2y_2 - \beta z_2 + \alpha_2(z_1 - z_2)
\end{cases}
\label{eq:mainsyst}
\end{equation}
where $(x_1, y_1, z_1)$ and $(x_2, y_2, z_2)$ denote state variables of the first and second Lorenz oscillators, respectively. The classical Lorenz parameters are fixed at $\sigma = 10$, $\rho = 28$, $\beta = 8/3$, ensuring chaotic behavior in the uncoupled case. The novel aspect of our model resides in the independent coupling coefficients $\alpha_1, \alpha_2 \in \mathbb{R}$, governing interaction strength and directionality between subsystems.
\\
The system manifests several important special cases:
\begin{itemize}
\item Decoupled systems: $\alpha_1 = \alpha_2 = 0$ yields two independent Lorenz oscillators.
\item {Symmetric coupling:} $\alpha_1 = \alpha_2 \neq 0$ produces reciprocal interactions.
\item {Asymmetric coupling:} $\alpha_1 \neq \alpha_2$ enables non-reciprocal interactions, representing realistic scenarios.
\end{itemize}
Finally, to demonstrate the practical relevance of the proposed coupled hyperchaotic system, a secure communication framework based on chaos synchronization is considered. The encryption/decryption example is not intended to introduce a new cryptographic architecture; rather, it serves to illustrate how the hyperchaotic dynamics and synchronization properties of the proposed system can be exploited in secure transmission applications.

The remainder of this manuscript is organized as follows: Section \ref{sec:model}  provides comprehensive dynamical analysis including dissipativity conditions, equilibrium points, local stability properties, Lyapunov spectra, and fractal dimensions. Section  \ref{sec3} develops theoretical conditions for complete synchronization between subsystems. Section \ref{sec4} presents numerical verification of synchronization phenomena through various parameter settings and simulations. Section \ref{sec5} provides critical discussion of results and comparative analysis of synchronization performance, and Section \ref{sec6} concludes the paper.

\section{Dynamical Analysis and Mathematical Formulation}\label{sec:model}
\subsection{Dissipativity Analysis}
A necessary condition for the existence of bounded chaotic attractors in
continuous-time dynamical systems is dissipativity, which ensures the
contraction of phase-space volumes as time evolves. For the proposed
six-dimensional coupled Lorenz system, the divergence of the vector
field $F$ is given by:
\begin{equation*}
\nabla \cdot F = -2\sigma - 2 - 2\beta - 3\alpha_{1} - 3\alpha_{2}.
\end{equation*}
With classical Lorenz parameters ($\sigma = 10,\ \beta = \frac{8}{3}$),
the system exhibits dissipative behavior when $\nabla \cdot F < 0$,
requiring:
\begin{equation*}
\alpha_{1} + \alpha_{2} > -\frac{82}{9}.
\end{equation*}
This condition is easily satisfied for practical coupling strengths.
Consequently, all trajectories are ultimately confined within a bounded
invariant set, ensuring the contraction of phase-space volumes and
providing a necessary foundation for the emergence of chaotic and
hyperchaotic dynamics.

\subsection{Equilibrium Points and Local Stability Analysis}

The equilibrium points are obtained by imposing
$\dot{x}_{1} = \dot{y}_{1} = \dot{z}_{1} = \dot{x}_{2} = \dot{y}_{2} = \dot{z}_{2} = 0$.
Although the proposed model involves asymmetric diffusive couplings
$(\alpha_{1}\neq\alpha_{2})$, the stationary solutions are found to be
synchronous and coincide with those of the classical Lorenz system.

\begin{figure}[H]
\centering
\begin{tabular}{cc}
    (c) & (d) \\
    \includegraphics[width=10cm,height=5.2cm]{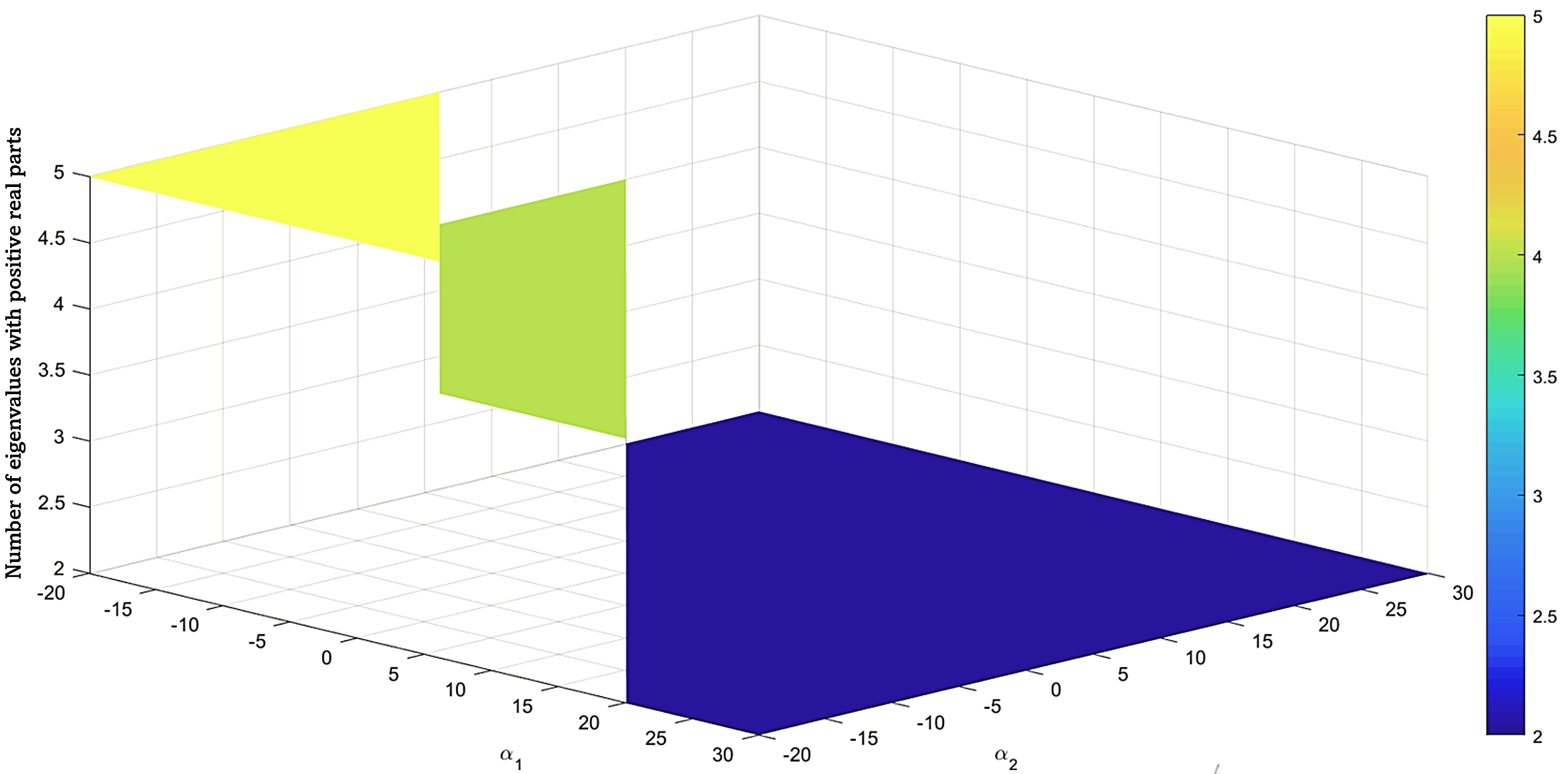} &
    \includegraphics[width=5.2cm,height=5.2cm]{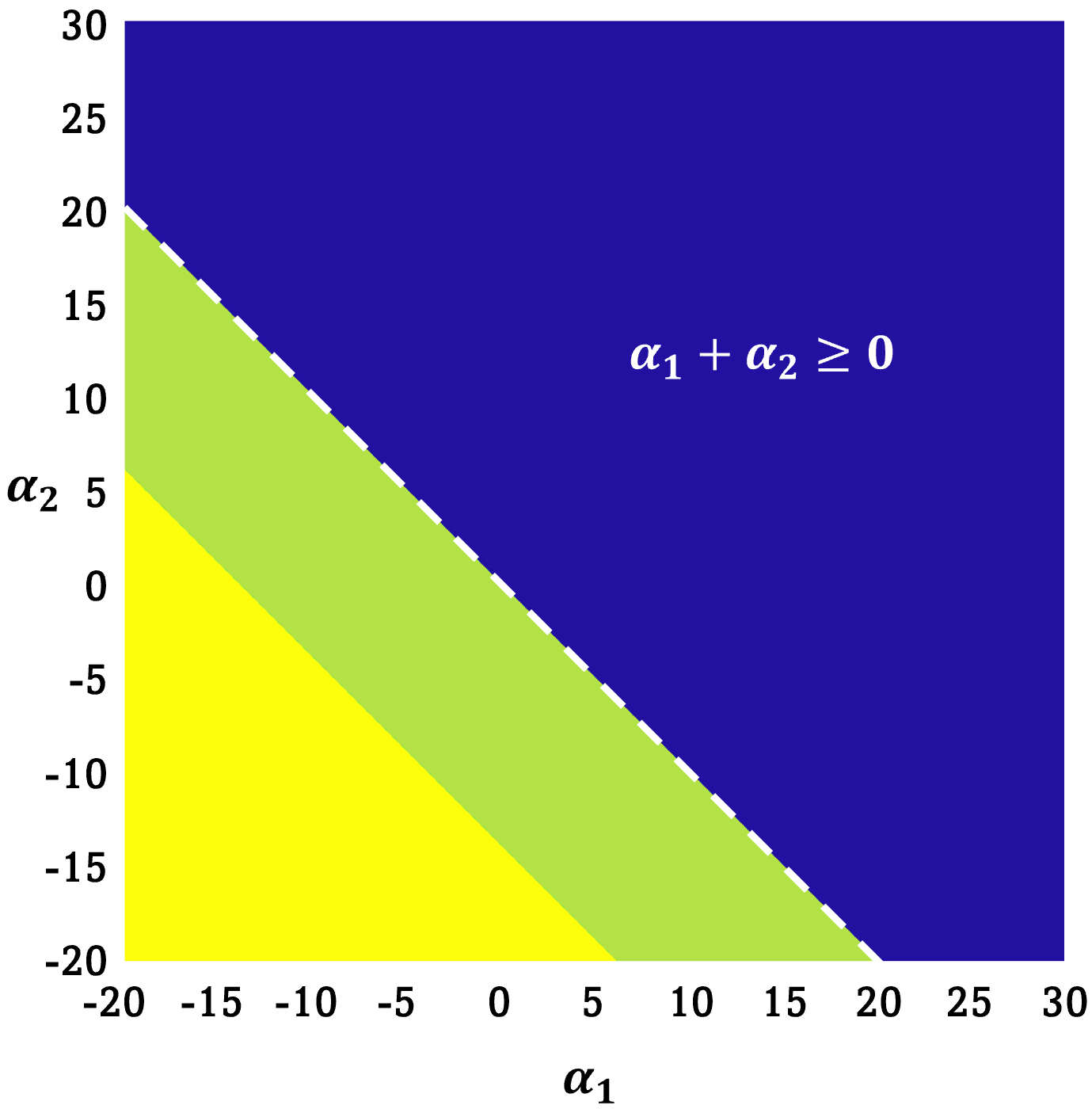}
\end{tabular}
\caption{Number of eigenvalues with positive real parts associated with the equilibrium points \( E_0 \) and \( E_+ \) as a function of the coupling coefficients \( \alpha_1 \) and \( \alpha_2 \). 3D representations are shown in (a) and (c), while the corresponding 2D maps are presented in (b) and (d), respectively.}
\label{X}
\end{figure}

The system admits the following equilibrium points:

\begin{itemize}
\item Trivial equilibrium: $E_{0} = (0,0,0,0,0,0)$. The Jacobian matrix
  evaluated at $E_{0}$ is given by:
\end{itemize}

\begin{equation*}
J_{\text{coupled}}(E_{0}) = \begin{pmatrix}
 -\sigma - \alpha_{1} & \sigma & 0 & \alpha_{1} & 0 & 0 \\
\rho & -1 - \alpha_{1} & 0 & 0 & \alpha_{1} & 0 \\
0 & 0 & -\beta - \alpha_{1} & 0 & 0 & \alpha_{1} \\
\alpha_{2} & 0 & 0 & -\sigma - \alpha_{2} & \sigma & 0 \\
0 & \alpha_{2} & 0 & \rho & -1 - \alpha_{2} & 0 \\
0 & 0 & \alpha_{2} & 0 & 0 & -\beta - \alpha_{2}
\end{pmatrix}.
\end{equation*}
\begin{itemize}
\item Non-trivial equilibria:
\end{itemize}
\begin{equation*}
E_{\pm} = \left( \pm \sqrt{\beta(\rho - 1)}, \pm \sqrt{\beta(\rho - 1)},\rho - 1, \pm \sqrt{\beta(\rho - 1)}, \pm \sqrt{\beta(\rho - 1)},\rho - 1 \right).
\end{equation*}
The corresponding Jacobian matrices are:
\begin{equation*}
J_{\text{coupled}}(E_{\pm}) = \begin{pmatrix}
 -\sigma - \alpha_{1} & \sigma & 0 & \alpha_{1} & 0 & 0 \\
1 & -1 - \alpha_{1} & \mp s & 0 & \alpha_{1} & 0 \\
\pm s & \pm s & -\beta - \alpha_{1} & 0 & 0 & \alpha_{1} \\
\alpha_{2} & 0 & 0 & -\sigma - \alpha_{2} & \sigma & 0 \\
0 & \alpha_{2} & 0 & 1 & -1 - \alpha_{2} & \mp s \\
0 & 0 & \alpha_{2} & \pm s & \pm s & -\beta - \alpha_{2}
\end{pmatrix},
\end{equation*}
where $s = \sqrt{\beta(\rho - 1)}$.
\\
For the classical Lorenz parameters $\beta = \frac{8}{3}$ and $\rho = 28$,
$s = \sqrt{\frac{8}{3}(27)} \approx 8.4853$.
\begin{equation*}
J(E_{+}) = \begin{pmatrix}
 -\sigma - \alpha_{1} & \sigma & 0 & \alpha_{1} & 0 & 0 \\
1 & -1 - \alpha_{1} & -8.4853 & 0 & \alpha_{1} & 0 \\
8.4853 & 8.4853 & -\beta - \alpha_{1} & 0 & 0 & \alpha_{1} \\
\alpha_{2} & 0 & 0 & -\sigma - \alpha_{2} & \sigma & 0 \\
0 & \alpha_{2} & 0 & 1 & -1 - \alpha_{2} & -8.4853 \\
0 & 0 & \alpha_{2} & 8.4853 & 8.4853 & -\beta - \alpha_{2}
\end{pmatrix},
\end{equation*}
\vspace{0.5cm}
\begin{equation*}
J(E_{-}) = \begin{pmatrix}
 -\sigma - \alpha_{1} & \sigma & 0 & \alpha_{1} & 0 & 0 \\
1 & -1 - \alpha_{1} & 8.4853 & 0 & \alpha_{1} & 0 \\
-8.4853 & -8.4853 & -\beta - \alpha_{1} & 0 & 0 & \alpha_{1} \\
\alpha_{2} & 0 & 0 & -\sigma - \alpha_{2} & \sigma & 0 \\
0 & \alpha_{2} & 0 & 1 & -1 - \alpha_{2} & 8.4853 \\
0 & 0 & \alpha_{2} & -8.4853 & -8.4853 & -\beta - \alpha_{2}
\end{pmatrix}.
\end{equation*}
The local stability of the equilibrium points is investigated through
the eigenvalues $\lambda_{i}$ of the Jacobian matrices, obtained from
$\det(J - \lambda I) = 0$. Because the resulting characteristic equation
is a sixth-order polynomial whose coefficients depend on two independent
coupling parameters, an analytical solution becomes extremely
cumbersome. Therefore, a numerical approach based on Matlab's \texttt{eig()}
function is employed to evaluate the eigenvalue spectrum over a wide
range of $(\alpha_{1}, \alpha_{2})$. The results are summarized in Figure \ref{X}, which presents the number of
eigenvalues with positive real parts for the equilibrium points $E_{0}$
and $E_{+}$. The results obtained for $E_{-}$ are similar to those of
$E_{+}$ and are therefore omitted.
\begin{figure}[H]
\centering
\begin{tabular}{cc}
    (a) & (b)  \\
    \includegraphics[width=10cm,height=5.2cm]{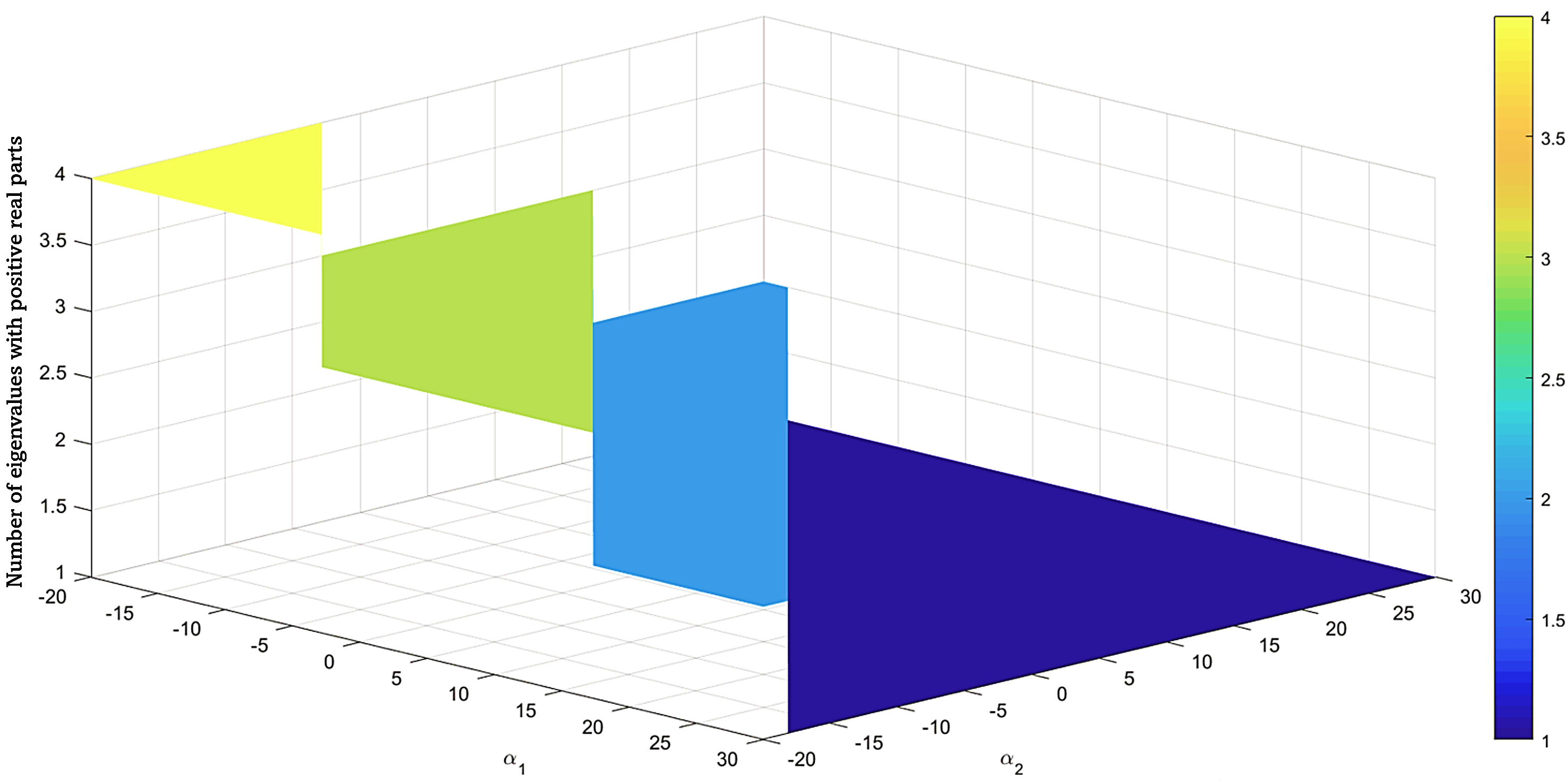} &
    \includegraphics[width=5.2cm,height=5.2cm]{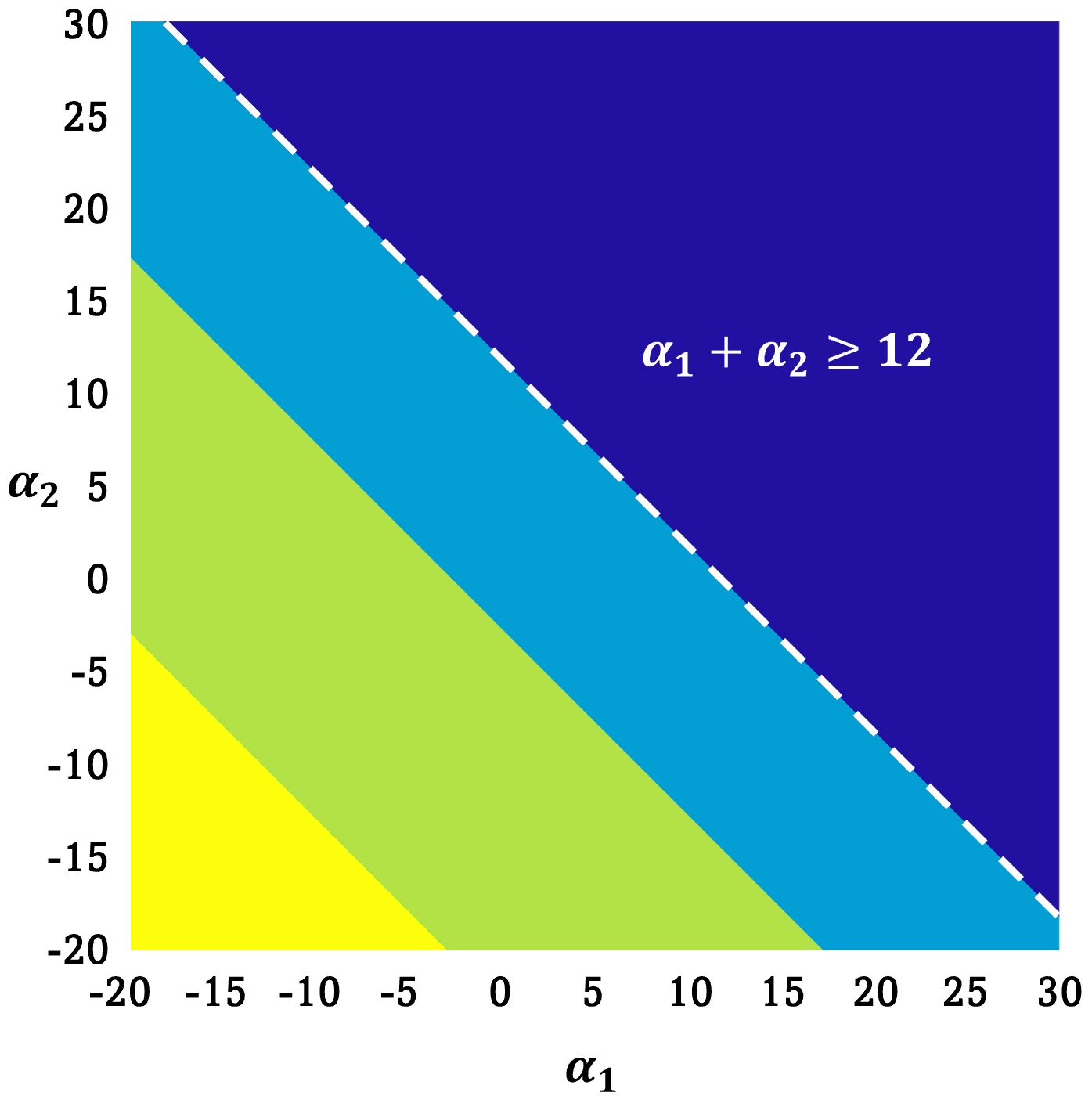} 
\end{tabular}
\end{figure}
\noindent
Figure \ref{X} reveals that negative coupling coefficients significantly
increase the number of eigenvalues with positive real parts, reaching up
to four unstable directions around the equilibrium points. This behavior
indicates stronger local instability and suggests the possibility of
richer dynamical regimes. As the coupling parameters increase toward
positive values, the number of positive eigenvalues progressively
decreases from four to one, demonstrating that positive diffusive
coupling reduces the number of locally unstable directions. These observations are consistent with the divergence of the proposed
system $\left( \nabla \cdot F = -\frac{82}{3} - 3(\alpha_1 + \alpha_2) \right)$.\\
Indeed, increasing positive coupling strengths makes the divergence more
negative, thereby enhancing phase-space contraction and overall
dissipation. Conversely, negative coupling strengths reduce the
contraction rate, allowing stronger local expansion mechanisms around
the equilibrium points.

It should be emphasized that this analysis only characterizes the local
stability properties near the equilibrium points. The actual global
dynamical behavior of the system is investigated in the next subsection
through the computation of Lyapunov exponents.

\subsection{Lyapunov Spectrum and Hyperchaotic Dynamics}

Hyperchaotic behavior necessitates at least two positive Lyapunov exponents (LEs). Lyapunov exponents provide a quantitative measure of the average
exponential divergence or convergence of neighboring trajectories in
phase space and constitute one of the most reliable indicators of
chaotic behavior.

For a nonlinear dynamical system $\dot{X} = F(X)$, an infinitesimal
perturbation evolves according to
$\| \delta X(t) \| \approx e^{\lambda t} \| \delta X(0) \|$.
The associated Lyapunov exponents are defined as:
\begin{equation*}
\lambda_{i} = \lim_{t \rightarrow \infty}\frac{1}{t}\ln\frac{\| \delta X_{i}(t) \|}{\| \delta X_{i}(0) \|}.
\end{equation*}
A positive largest Lyapunov exponent indicates chaotic behavior, whereas
hyperchaotic dynamics require the existence of at least two positive
Lyapunov exponents.
\\
The Lyapunov spectrum was computed numerically for varying values of the
coupling coefficients $(\alpha_{1}, \alpha_{2})$. For each parameter
pair, the coupled Lorenz equations and their associated variational
equations were integrated simultaneously using a fourth-order
Runge--Kutta (RK4) scheme. The RK4 method was selected due to its
superior numerical accuracy and stability when dealing with highly
sensitive chaotic trajectories.
\begin{figure}[H]
\centering
\begin{tabular}{cc}
    (a) & (b)  \\
    \includegraphics[width=10.7cm,height=5.4cm]{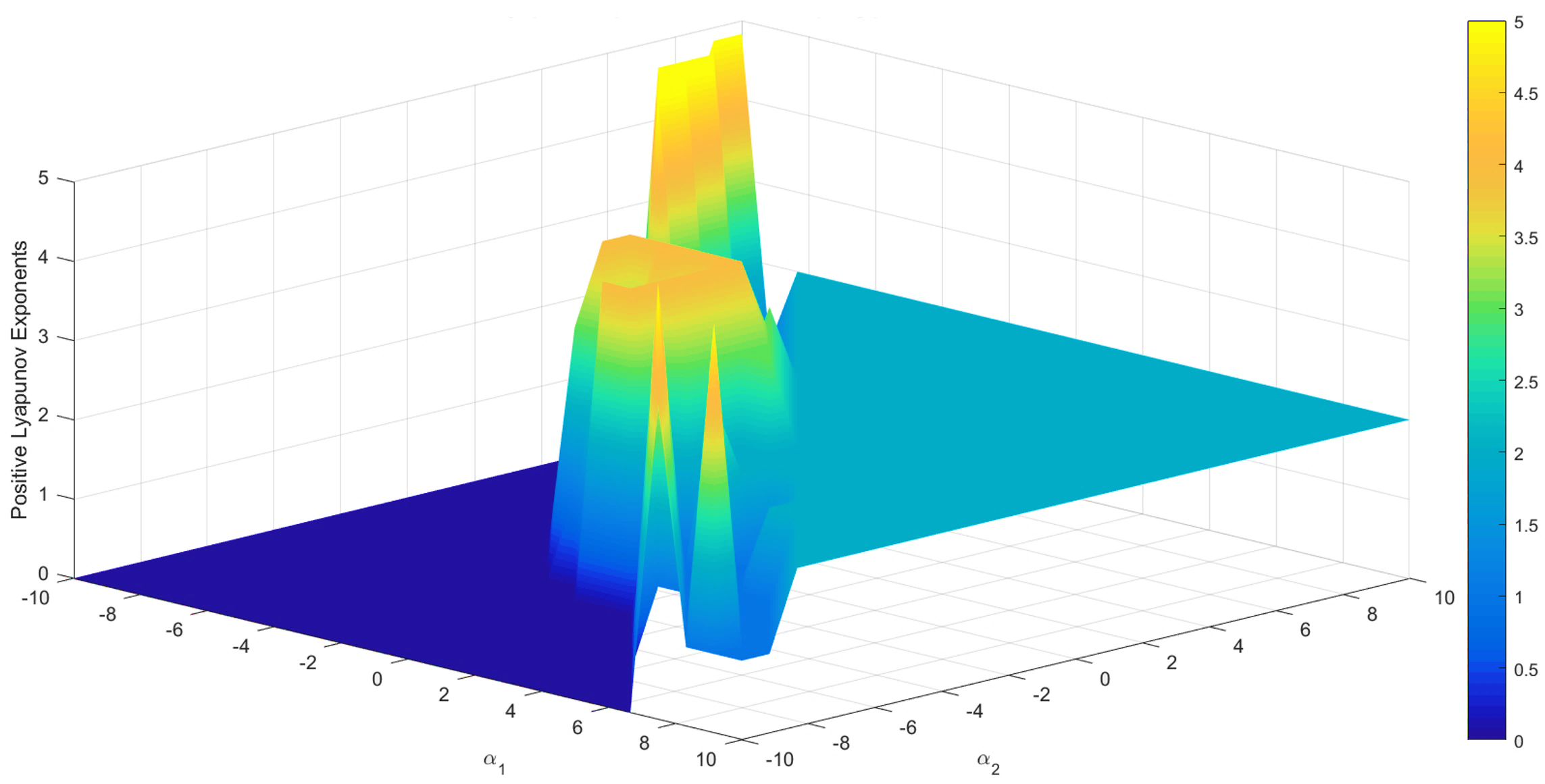} &
    \includegraphics[width=5.2cm,height=5.2cm]{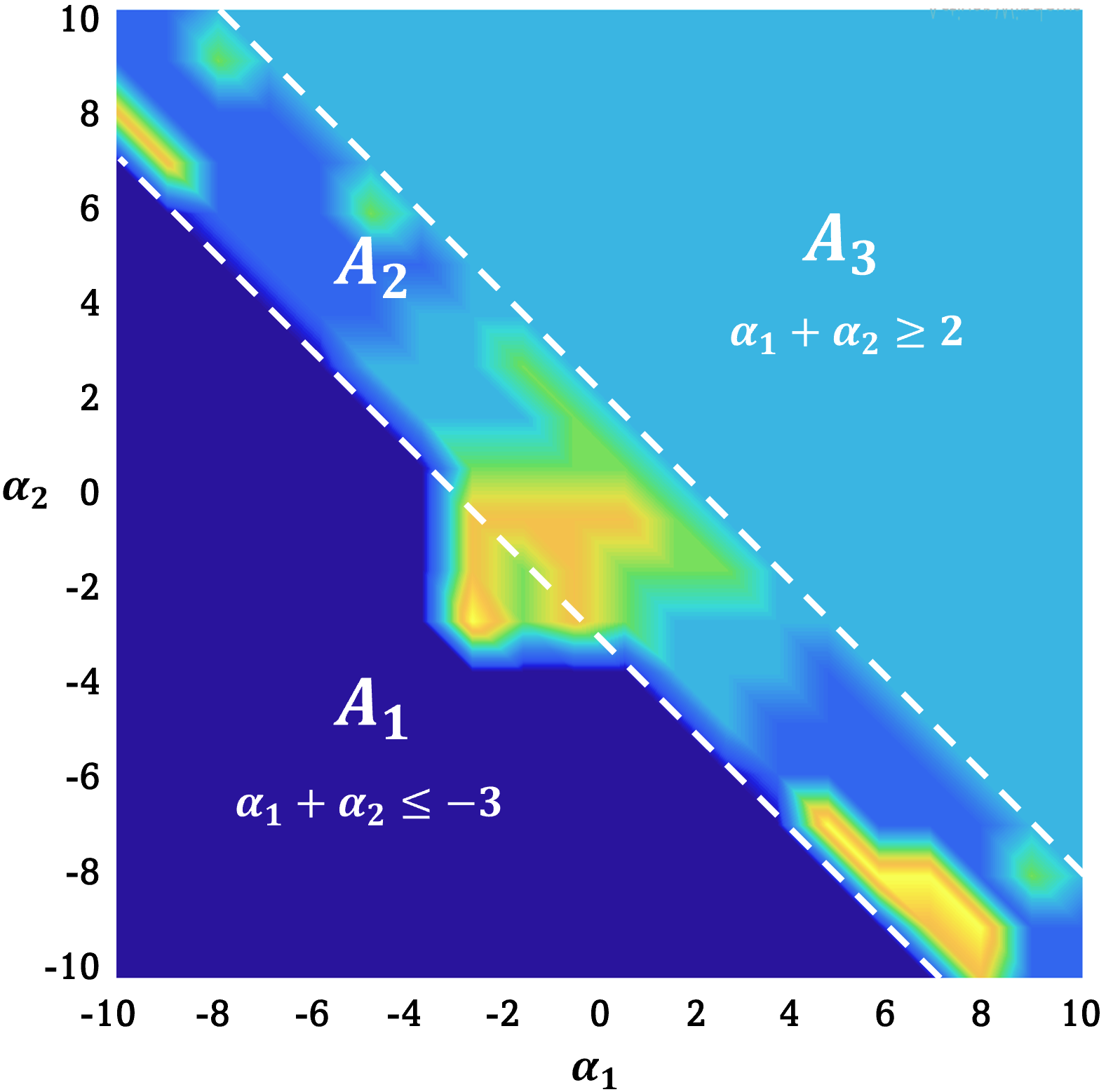} 
\end{tabular}
\caption{Lyapunov exponent spectrum: (a) 3D view, (b) 2D view}\label{Y}
\end{figure}
\noindent
As illustrated in Figure \ref{Y} three distinct dynamical regions can be
identified. Region $(\mathbf{A_1})$, corresponding to $\alpha_{1} + \alpha_{2} \leq -3$,
is characterized by Lyapunov spectra that are predominantly negative, except
for a small subset of parameter combinations located near the center of the
region. This indicates the disappearance of sustained chaotic activity and
the emergence of stable non-chaotic dynamics. Region $(\mathbf{A_2})$ constitutes a
transition zone in which the number of positive Lyapunov exponents varies
significantly with the coupling parameters. Consequently, the system
alternates between chaotic and hyperchaotic behaviors. Finally, Region $(\mathbf{A_3})$,
defined by $\alpha_{1} + \alpha_{2} \geq 2$, exhibits two positive Lyapunov
exponents throughout the investigated domain. This confirms the persistence
of robust hyperchaotic dynamics under sufficiently strong positive coupling.

The Lyapunov analysis therefore complements the local stability study by
providing a global characterization of the system dynamics and revealing
the parameter regions associated with non-chaotic, chaotic, and
hyperchaotic behaviors.

\subsection{Bifurcation analysis}
\begin{figure}[H]
\centering
\begin{tabular}{cc}
    (a) & (b) \\
    \includegraphics[width=10.7cm,height=5.9cm]{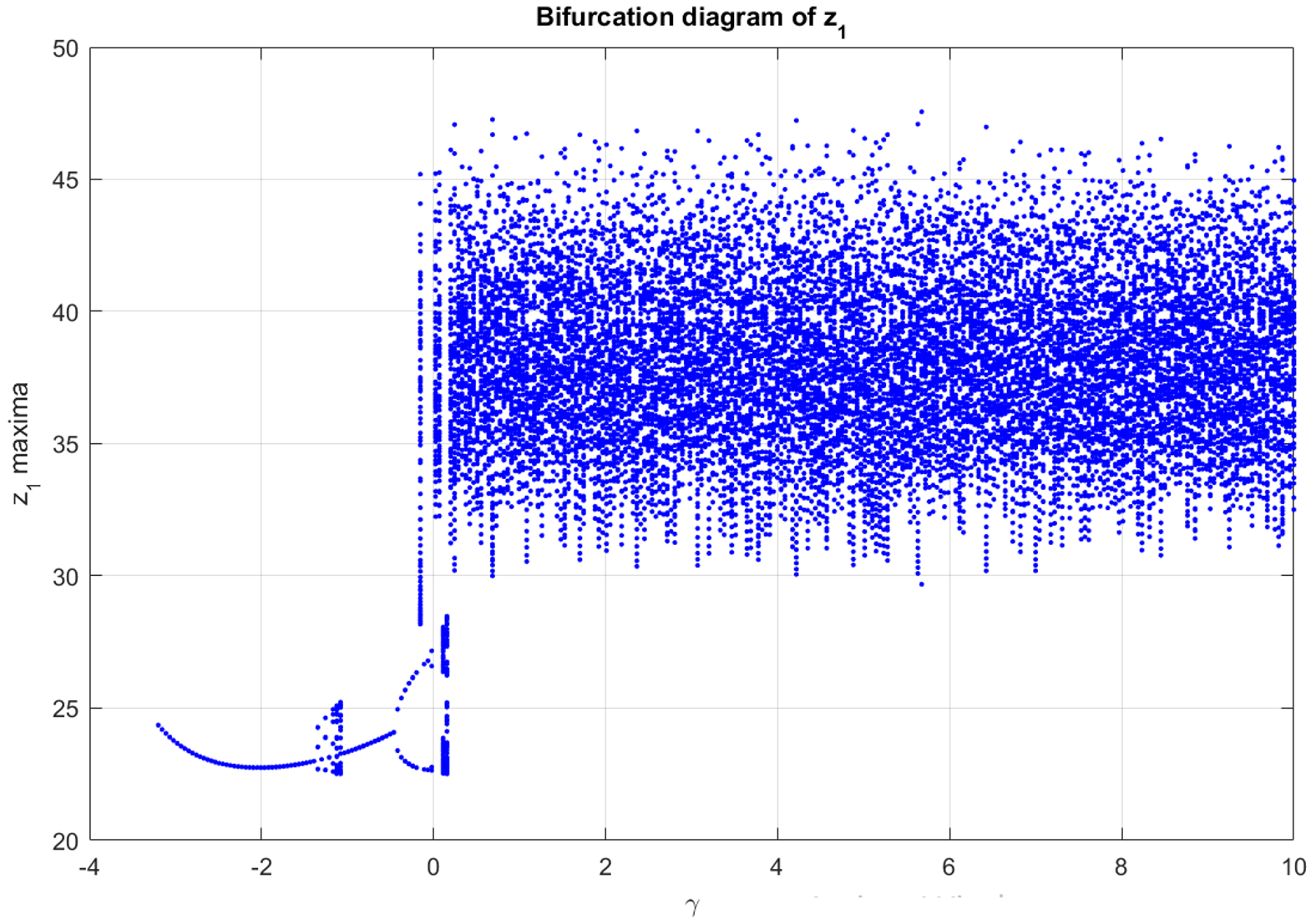} &
    \includegraphics[width=5.2cm,height=5.7cm]{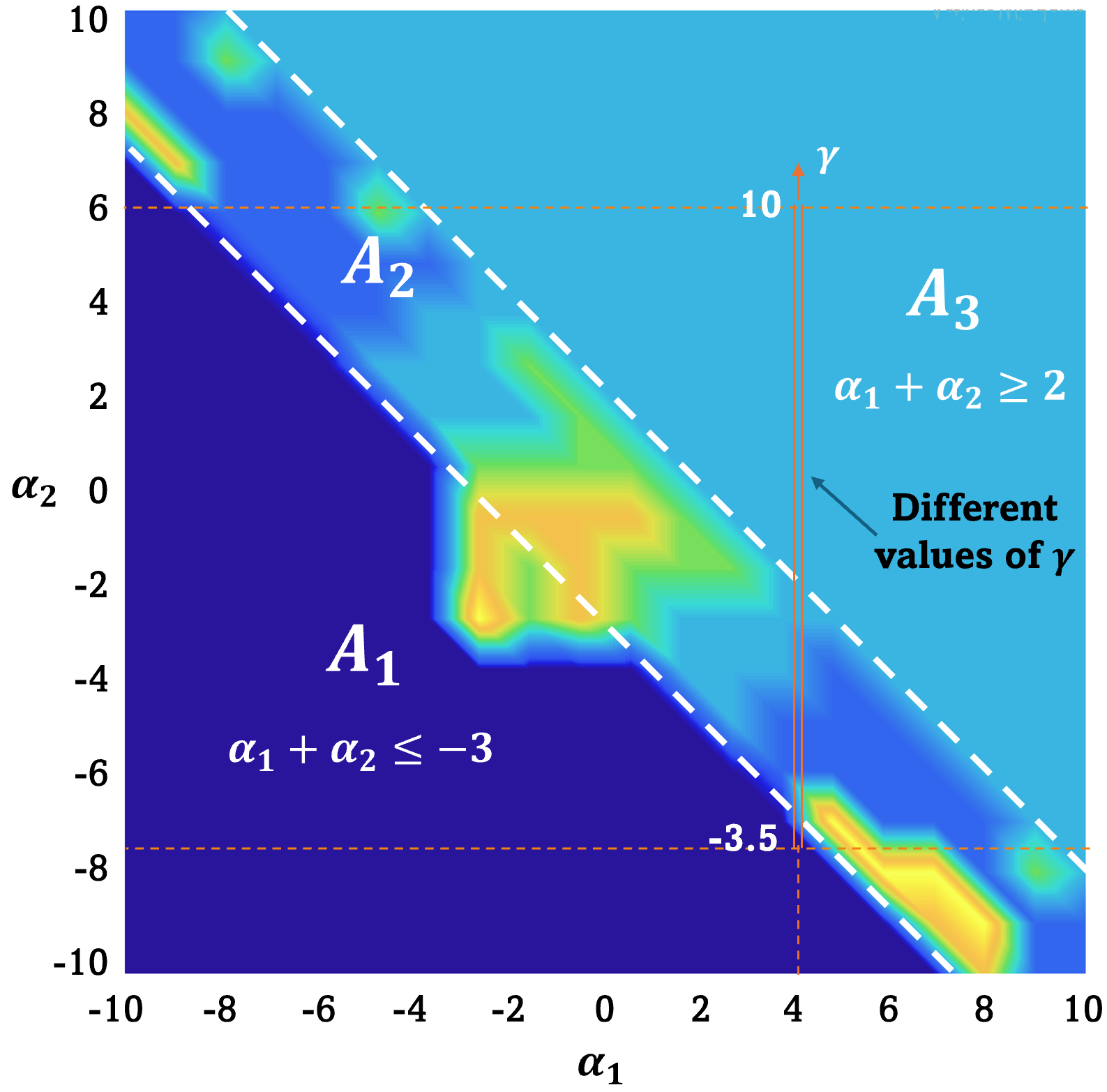} 
\end{tabular}
\caption{(a) Bifurcation diagram of $z_1$ as a function of the coupling parameter $\gamma$. (b) Identification of the investigated coupling parameter range ($-3.5 \leq \gamma \leq 10$) within the dynamical regions of the proposed system.}
\label{fig:bifurcation}
\end{figure}

To further validate the parameter regions identified by the Lyapunov exponent analysis, the bifurcation diagram of the state variable $z_1$ was generated by varying the overall coupling strength $\gamma = \alpha_1 + \alpha_2$, while fixing $\alpha_1 = 4$ and varying $\alpha_2$, as shown in Fig.~\ref{fig:bifurcation}. The variable $z_1$ was selected because it provides a representative description of the global dynamics of the proposed six-dimensional system. Similar bifurcation patterns were observed for the remaining state variables.

The obtained bifurcation diagram confirms the dynamical transitions predicted by the Lyapunov exponent analysis. For negative coupling strengths, successive bifurcations are observed, revealing a gradual transition from predominantly non-chaotic dynamics toward chaotic behavior. This result is consistent with Regions $(\mathbf{A_1})$ and $(\mathbf{A_2})$  of the Lyapunov spectrum, where most Lyapunov exponents are negative, except for a narrow subset of parameter values exhibiting weak chaotic dynamics. As $\gamma$ increases, the number and density of extrema progressively increase, indicating the emergence of fully developed chaotic oscillations. For positive coupling strengths, the bifurcation diagram evolves into a broad continuous chaotic branch, demonstrating that the proposed asymmetrically coupled Lorenz system preserves robust chaotic dynamics over a wide interval of coupling strengths.
\section{Synchronization Analysis}\label{sec3}
Synchronization of chaotic oscillators has attracted considerable
attention due to its numerous applications in secure communications,
neural networks, biological systems, and information processing. In this
section, sufficient conditions ensuring complete synchronization of the
proposed coupled Lorenz system are established using the Lyapunov direct
method.

\subsection{Error Dynamical System}

To investigate synchronization between the two coupled oscillators, the
synchronization error variables are defined as
$e_{1} = x_{2} - x_{1},\; e_{2} = y_{2} - y_{1},\; e_{3} = z_{2} - z_{1}$.
Complete synchronization is achieved when
$\lim_{t \rightarrow \infty} e_{i}(t) = 0$, for $i = 1,2,3$.
Therefore, the synchronization problem reduces to analyzing the
asymptotic stability of the origin $(e = 0)$ of the error dynamical
system. By subtracting the state equations of the first subsystem from
those of the second subsystem and rearranging the resulting terms, the
error dynamics can be expressed as:
\begin{equation*}
\left\{
\begin{array}{l}
\dot{e}_{1} = -(\sigma + \alpha_{1} + \alpha_{2})e_{1} + \sigma e_{2} \\[2mm]
\dot{e}_{2} = (\rho - z_{2})e_{1} - (1 + \alpha_{1} + \alpha_{2})e_{2} - x_{1}e_{3} \\[2mm]
\dot{e}_{3} = y_{2}e_{1} + x_{1}e_{2} - (\beta + \alpha_{1} + \alpha_{2})e_{3}
\end{array}
\right.
\end{equation*}
Introducing the parameter $\gamma = \alpha_{1} + \alpha_{2}$, the error
system can be rewritten in compact matrix form as:
\begin{equation*}
\begin{bmatrix}
\dot{e}_{1} \\ \dot{e}_{2} \\ \dot{e}_{3}
\end{bmatrix}
=
\begin{bmatrix}
-(\sigma + \gamma) & \sigma & 0 \\
\rho - z_{2} & -(1 + \gamma) & -x_{1} \\
y_{2} & x_{1} & -(\beta + \gamma)
\end{bmatrix}
\begin{bmatrix}
e_{1} \\ e_{2} \\ e_{3}
\end{bmatrix}
= Q e
\end{equation*}
where $e = (e_{1},e_{2},e_{3})^{T}$ denotes the synchronization error
vector and $Q$ is given by
\begin{equation*}
Q = \begin{bmatrix}
-(\sigma + \gamma) & \sigma & 0 \\
\rho - z_{2} & -(1 + \gamma) & -x_{1} \\
y_{2} & x_{1} & -(\beta + \gamma)
\end{bmatrix}
\end{equation*}
\subsection{Lyapunov-Based Synchronization Conditions}
Synchronization of the proposed coupled system is achieved when the
synchronization errors asymptotically vanish. To establish this
property, consider the Lyapunov candidate function
\begin{equation*}
V(e) = \frac{1}{2}e^{T} P e,\quad V(e) > 0\ \forall e \neq 0,
\end{equation*}
where $P = P^{T} > 0$ is a symmetric positive-definite matrix.
\\
The time derivative of $V(e)$ along the trajectories of the error system
is given by
\begin{equation*}
\dot{V}(e) = e^{T} P \dot{e} = e^{T} P Q e.
\end{equation*}
The asymptotic convergence of the synchronization error is guaranteed if
$e^{T} P Q e < 0$.
\\
By applying Sylvester's criterion, the synchronization conditions of
the proposed coupled system are determined from the negative
definiteness of the matrix $(PQ)$, i.e., $PQ < 0$.
\\
Since the matrix $Q$ is generally non-symmetric, the previous condition
is replaced by the equivalent symmetric condition
$Q^{T} P + P Q < 0$.
\\
By choosing $P = \operatorname{diag}(p_{1},p_{2},p_{3})$ with
$p_{1},p_{2},p_{3} > 0$, the negative definiteness of $Q^{T}P + PQ$ is
equivalent to the alternating-sign conditions: $\Delta_{1} < 0$,
$\Delta_{2} > 0$ and $\Delta_{3} < 0$,
\\
where
\begin{equation*}
\begin{aligned}
\Delta_1 &= -2p_1(\sigma + \gamma), \\
\Delta_2 &= \det \begin{pmatrix}
-2p_1(\sigma + \gamma) & p_1\sigma + p_2(\rho - z_2), \\
p_1\sigma + p_2(\rho - z_2) & -2p_2(1 + \gamma)
\end{pmatrix} \\
\Delta_3 &= \det \begin{pmatrix}
-2p_1(\sigma + \gamma) & p_1\sigma + p_2(\rho - z_2) & p_3y_2 \\
p_1\sigma + p_2(\rho - z_2) & -2p_2(1 + \gamma) & (p_3 - p_2)x_1 \\
p_3y_2 & (p_3 - p_2)x_1 & -2p_3(\beta + \gamma)
\end{pmatrix}.
\end{aligned}
\end{equation*}
We obtain,
\begin{equation*}
\begin{aligned}
\Delta_1 &= -2p_1(\sigma + \gamma) \\
\Delta_2 &= 4p_1p_2(\sigma + \gamma)(1 + \gamma) - \bigl[p_1\sigma + p_2(\rho - z_2)\bigr]^2 \\
\Delta_3 &= -8p_1p_2p_3(\sigma + \gamma)(1 + \gamma)(\beta + \gamma) + 2p_1(\sigma + \gamma)(p_3 - p_2)^2 x_1^2 \\
&\qquad + 2p_3(\beta + \gamma)\bigl[p_1\sigma + p_2(\rho - z_2)\bigr]^2 + 2(p_3 - p_2)p_3x_1y_2\bigl[p_1\sigma + p_2(\rho - z_2)\bigr] + 2p_2(1 + \gamma)p_3^2 y_2^2.
\end{aligned}
\end{equation*}
If we choose $P = I$,  then $\Delta_1, \Delta_2$   and $\Delta_3$  simplify greatly:
\begin{equation*}
\begin{aligned}
\Delta_1 &= -2(\sigma + \gamma), \\
\Delta_2 &= 4\frac{\sigma + \gamma}{1 + \gamma} - \bigl[\sigma + (\rho - z_2)\bigr]^2, \\
\Delta_3 &= 2\left[-4(\sigma + \gamma)(1 + \gamma)(\beta + \gamma) + (\beta + \gamma)(\sigma + \rho - z_2)^2 + (1 + \gamma)y_2^2\right].
\end{aligned}
\end{equation*}
Since the state variables remain confined within bounded intervals,
i.e., $x_{\min} \leq x_{i}(t) \leq x_{\max}$,
$y_{\min} \leq y_{i}(t) \leq y_{\max}$,
$z_{\min} \leq z_{i}(t) \leq z_{\max}$, $i = 1,2$,
the stability conditions can be investigated over the effective ranges
of the trajectories rather than by using global worst-case bounds. Such
an approach yields less conservative synchronization criteria while
preserving the validity of the theoretical analysis.

\begin{figure}[H]
\centering
\includegraphics[width=14cm,height=7cm]{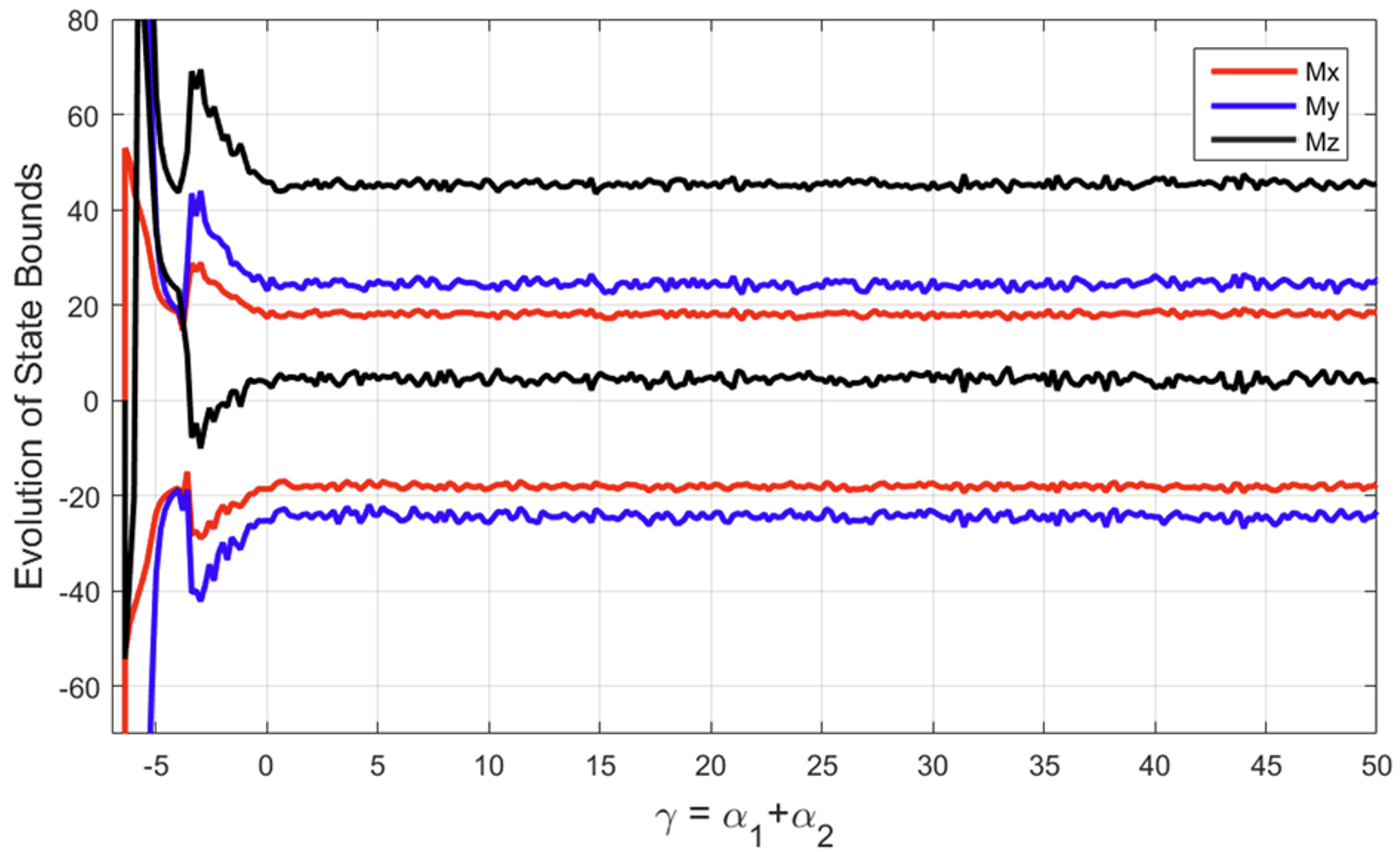}
\caption{Evolution of the state bounds as a function of the coupling
parameter $\gamma = \alpha_{1} + \alpha_{2}$. The upper and lower bounds
of the state variables $x$, $y$, and $z$ are extracted from the numerical
attractors for each value of $\gamma$.}\label{figZ}
\end{figure}

Based on the numerical analysis presented in Figure \ref{figZ}, the state variable
$z_{2}$ remains bounded within the interval $[0,48]$, whereas $y_{2}$
evolves within the range $[-28,28]$ over the investigated values of the
coupling parameter $\gamma$. These attractor-based bounds are subsequently
employed to derive less conservative sufficient conditions for $\Delta_{2}$
and $\Delta_{3}$.
Figure \ref{figZ} illustrates the evolution of the state bounds as a function of
the coupling parameter $\gamma = \alpha_{1} + \alpha_{2}$. For positive
values of $\gamma$, the attractor remains confined within nearly
constant bounds, indicating that the additional coupling mainly enhances
the dissipative nature of the system. Conversely, when $\gamma$ becomes
negative, the state amplitudes increase significantly due to the
reduction of the overall phase-space contraction rate. This behavior is
consistent with the divergence $\nabla \cdot F = -\frac{82}{3} - 3\gamma$,
whose magnitude decreases as $\gamma$ becomes more negative. Nevertheless,
the divergence remains negative over the considered range, ensuring the
persistence of a bounded attractor and preventing unbounded growth of the
trajectories. Consequently, negative coupling strengths enlarge the
attractor while preserving the dissipative character of the coupled
Lorenz system.

Based on the numerical bounds extracted from Figure \ref{figZ}, namely
$0 \leq z_{2}(t) \leq 48$, $-28 \leq y_{2}(t) \leq 28$ and
$\sigma = 10$, $\rho = 28$, we obtain
$\displaystyle\max_{z_{2} \in [0,48]} (\sigma + \rho - z_{2})^{2} = 38^{2} = 1444$.
\\
Consequently, a sufficient condition ensuring $\Delta_{2} > 0$ is
\begin{equation*}
4(10 + \gamma)(1 + \gamma) - 1444 > 0.
\end{equation*}
Solving the above inequality yields two roots. Since the negative root
does not belong to the hyperchaotic region identified in Figure \ref{Y}, only
the positive solution is retained. Therefore,
\begin{equation*}
\gamma > 14.03
\end{equation*}
To ensure that $\Delta_{3} < 0$, it is sufficient to impose
$(\beta + \gamma)(38 - z_{2})^{2} + (1 + \gamma)y_{2}^{2} < 4(\sigma + \gamma)(1 + \gamma)(\beta + \gamma)$.
\\
Using the numerical bounds extracted from Figure \ref{figZ},
$\max(38 - z_{2})^{2} = 1444$ and $\max y_{2}^{2} = 28^{2} = 784$ and
substituting the classical Lorenz parameters $\sigma = 10$, $\beta = \frac{8}{3}$,
the previous inequality becomes
\begin{equation*}
1444\left(\gamma + \frac{8}{3}\right) + 784(1 + \gamma) < 4(10 + \gamma)(1 + \gamma)\left(\gamma + \frac{8}{3}\right).
\end{equation*}
After simplification, a sufficient condition for $\Delta_{3} < 0$ is
obtained as
\begin{equation*}
4\gamma^{3} + 54.67\gamma^{2} - 2070.67\gamma - 2114.67 > 0.
\end{equation*}
Solving the above inequality yields a positive root approximately equal to:
\begin{equation*}
\gamma \approx 15.1
\end{equation*}
Finally, by jointly considering the dissipativity requirement, the
equilibrium-point stability analysis, the hyperchaotic operating region
identified through the Lyapunov spectrum, and the synchronization
conditions derived from the Lyapunov approach, a sufficient
synchronization criterion for the proposed coupled Lorenz system is
obtained as  $\alpha_{1} + \alpha_{2} > 15.1.$  More precisely, the trajectories of the two systems converge to each other:
\[
\lim_{t \to \infty} |x_1(t) - x_2(t)| = 0,
\qquad
\lim_{t \to \infty} |y_1(t) - y_2(t)| = 0,
\qquad
\lim_{t \to \infty} |z_1(t) - z_2(t)| = 0.
\]
It should be emphasized that the condition $\alpha_1+\alpha_2>15.1$, constitutes a sufficient synchronization condition derived from the Lyapunov framework and Sylvester's criterion. Therefore, it guarantees the asymptotic convergence of the synchronization error but does not represent a necessary condition. Consequently, synchronization may still occur for smaller values of the coupling parameter. This behavior is a well-known consequence of the conservatism inherent to Lyapunov-based sufficient criteria.

To illustrate this point and assess the conservatism of the theoretical bound, additional numerical investigations are presented as below. Several representative cases satisfying $\alpha_1+\alpha_2\leq 15.1$, are examined, showing that complete synchronization can still be achieved despite the fact that the sufficient analytical conditions are not fulfilled.
\section{Case of Fixed Parameters}\label{sec4}
\subsection{Synchronization case with $\alpha_1 = 1.5$ and $\alpha_2 = 4.5$}

\begin{figure}[H]
\centering
\includegraphics[width=14cm,height=6cm]{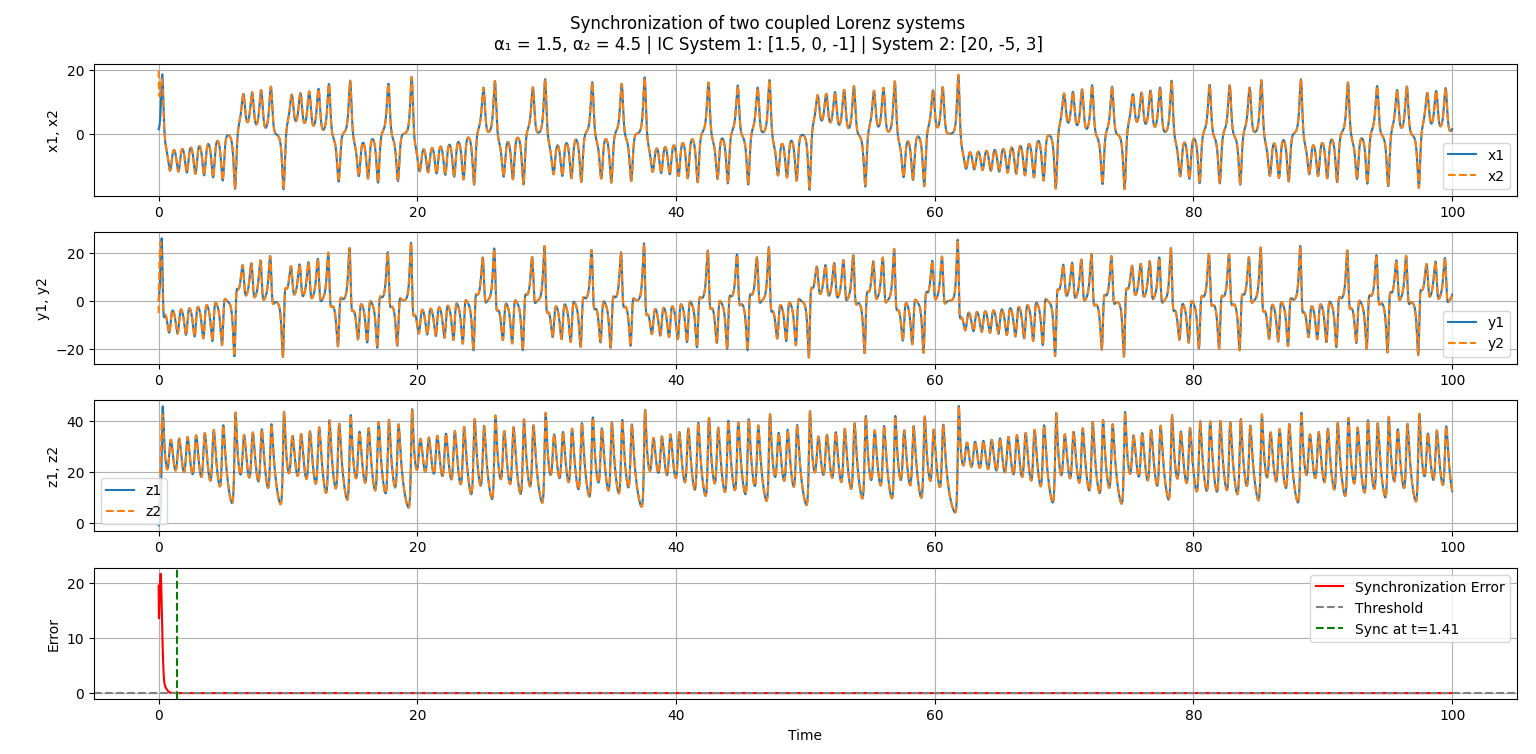}
\caption{Time series synchronization for $\alpha_1 = 1.5$, $\alpha_2 = 4.5$.}
\label{fig:di1}
\end{figure}

\begin{figure}[H]
\centering
\includegraphics[width=14cm,height=6cm]{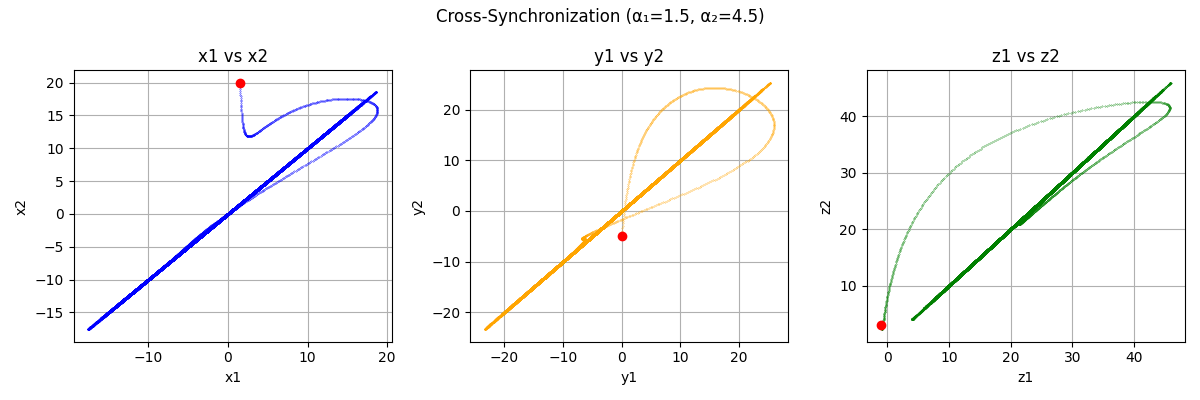}
\caption{Diagonal synchronization plot for $\alpha_1 = 1.5$, $\alpha_2 = 4.5$.}
\label{fig:di2}
\end{figure}

\begin{figure}[H]
\centering
\includegraphics[width=14cm,height=6cm]{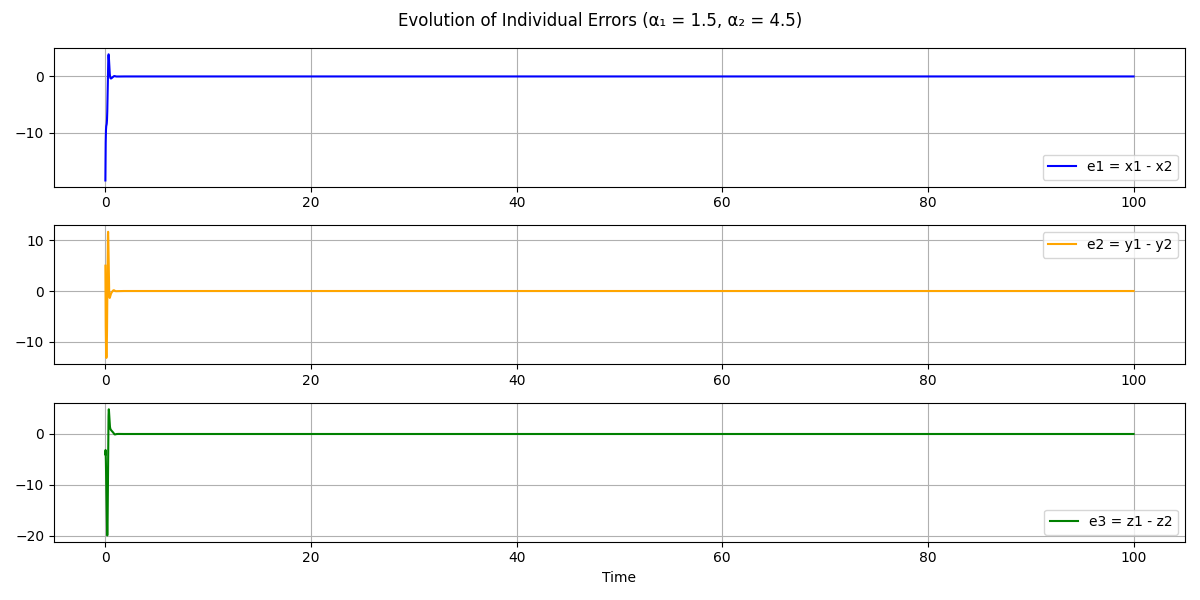}
\caption{Error convergence for $\alpha_1 = 1.5$, $\alpha_2 = 4.5$.}
\label{fig:di3}
\end{figure}

\begin{figure}[H]
\centering
\includegraphics[width=14cm,height=6cm]{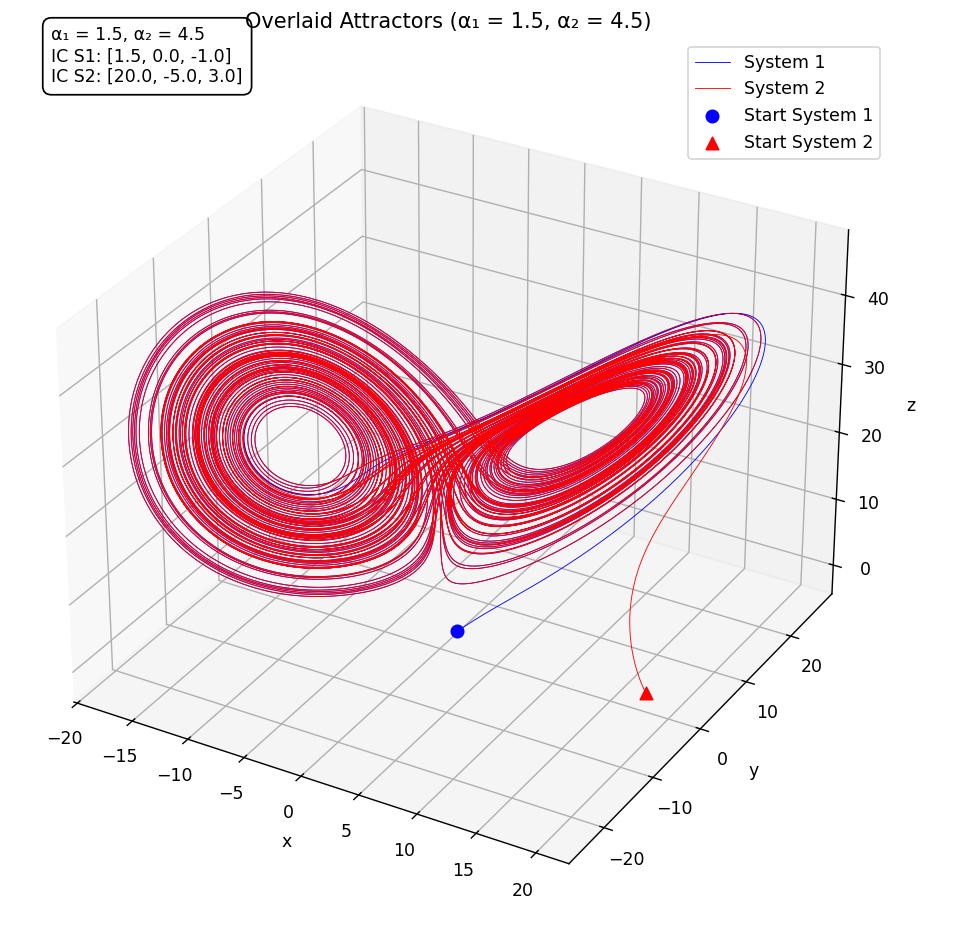}
\caption{Synchronized attractor for $\alpha_1 = 1.5$, $\alpha_2 = 4.5$.}
\label{fig:di5}
\end{figure}

\subsection{Non-Synchronization case with $\alpha_1 = 0.5$, $\alpha_2 = 0.2$}

\begin{figure}[H]
\centering
\includegraphics[width=14cm,height=5cm]{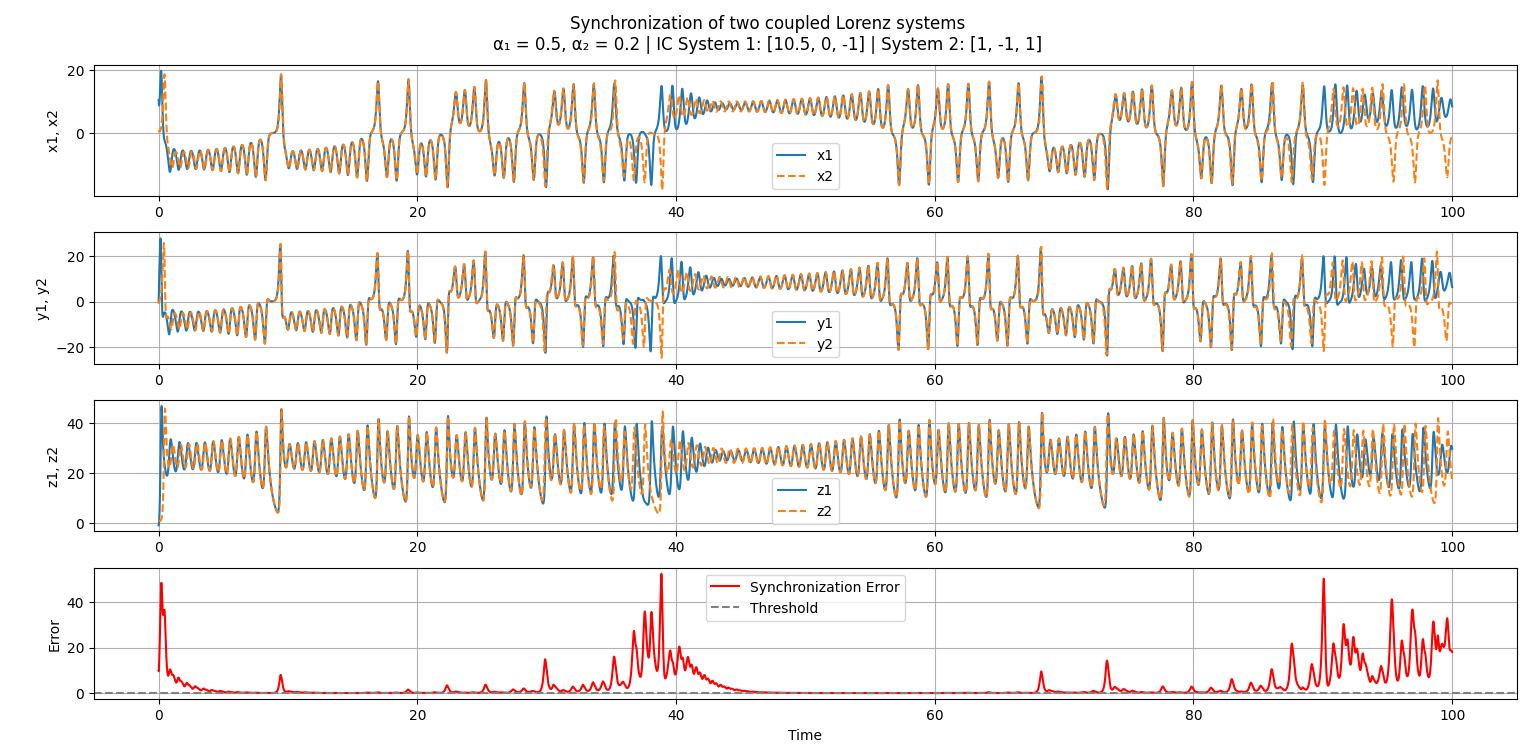}
\caption{Time series $ (x_1(t), x_2(t)) $, $ (y_1(t), y_2(t)) $, and $ (z_1(t), z_2(t)) $ of system \eqref{eq:mainsyst} showing non-convergence for $\alpha_1 = 0.5$, $\alpha_2 = 0.2$.}
\label{fig:Fi1}
\end{figure}

\begin{figure}[H]
\centering
\includegraphics[width=14cm,height=5cm]{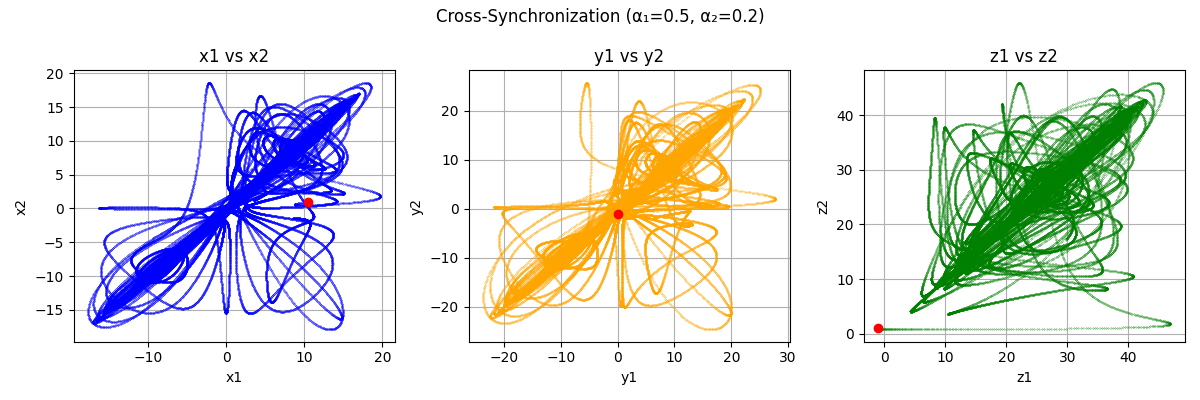}
\caption{Phase portrait showing non-synchronization for $\alpha_1 = 0.5$, $\alpha_2 = 0.2$.}
\label{fig:Fi2}
\end{figure}

\begin{figure}[H]
\centering
\includegraphics[width=14cm,height=5cm]{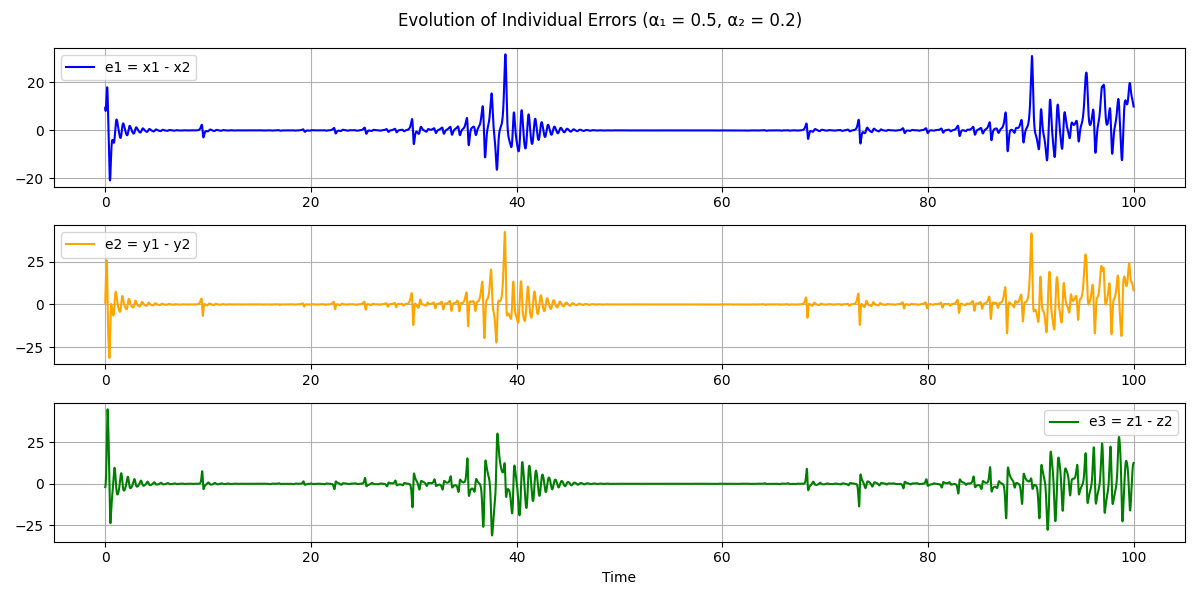}
\caption{Error evolution $e_1, e_2, e_3$ for non-synchronized case $\alpha_1 = 0.5$, $\alpha_2 = 0.2$.}
\label{fig:Fi3}
\end{figure}

\begin{figure}[H]
\centering
\includegraphics[width=12cm,height=5cm]{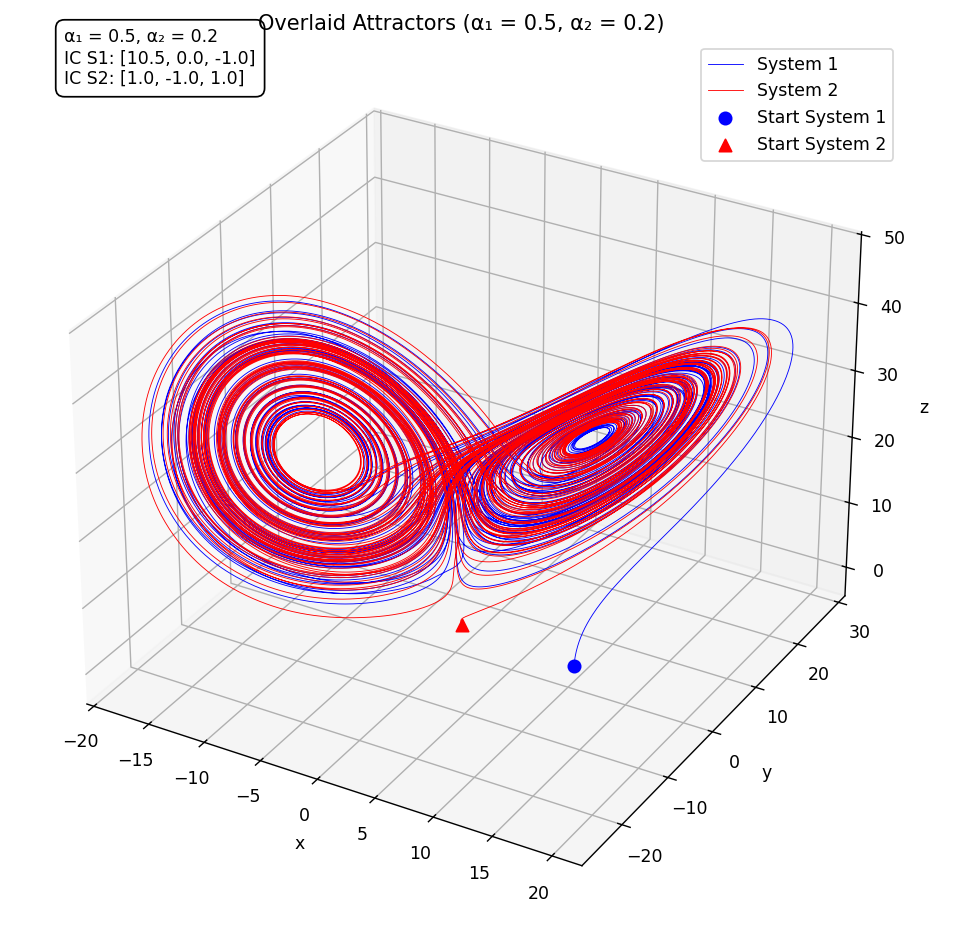}
\caption{Attractor showing non-synchronized behavior for $\alpha_1 = 0.5$, $\alpha_2 = 0.2$.}
\label{fig:Fi5}
\end{figure}

\subsection{Coupled Lorenz system with $\alpha_1 = 0.8$, $\alpha_2 = 0.2$}

\subsubsection{Lyapunov Exponent Spectrum and Hyperchaotic Behavior}
\begin{figure}[!ht]
\centering
\includegraphics[width=12cm,height=5cm]{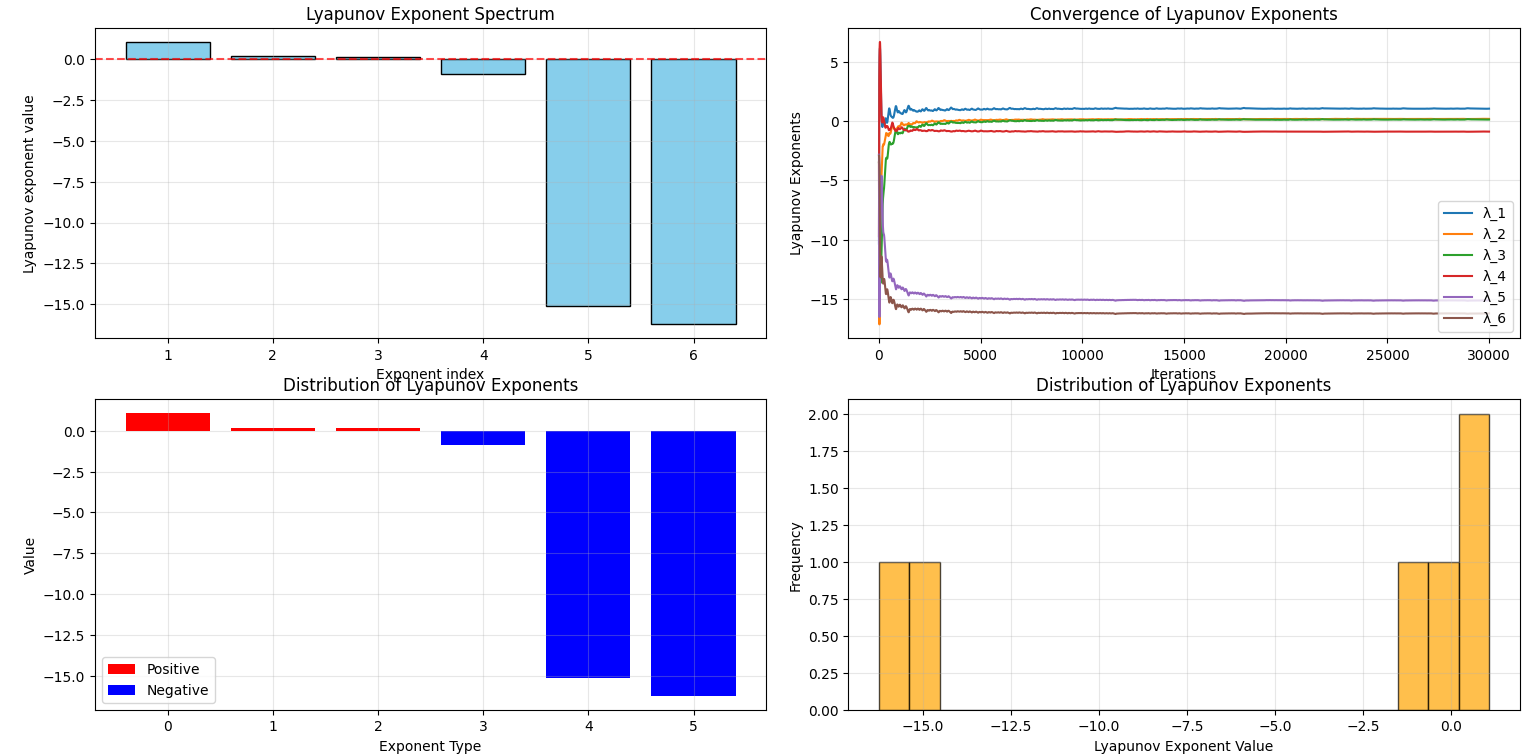}
\caption{Lyapunov exponent spectrum for coupled Lorenz system with $\alpha_1 = 0.8$, $\alpha_2 = 0.2$, demonstrating three positive exponents confirming hyperchaotic behavior.}
\label{fig:lyapunov_spectrum}
\end{figure}

Lyapunov exponents (LEs) quantify the average exponential rates of divergence or convergence of adjacent trajectories in phase space. Hyperchaotic behavior necessitates at least two positive LEs.
\\
As shown in Figure \eqref{fig:lyapunov_spectrum}, the Lyapunov exponent spectrum for coupled Lorenz system with $\alpha_1 = 0.8$, $\alpha_2 = 0.2$ demonstrates three positive exponents confirming hyperchaotic behavior.
We computed LEs using the Benettin algorithm over 30,000 iterations. For representative parameters $\alpha_1 = 0.8$, $\alpha_2 = 0.2$, we obtained:
\begin{equation}
\begin{aligned}
\lambda_1 &= 1.0679, \quad \lambda_2 = 0.2036, \quad \lambda_3 = 0.1479, \\
\lambda_4 &= -0.8717, \quad \lambda_5 = -15.1399, \quad \lambda_6 = -16.2390
\end{aligned}
\end{equation}

\begin{figure}[H]
\centering
\includegraphics[width=0.75\linewidth]{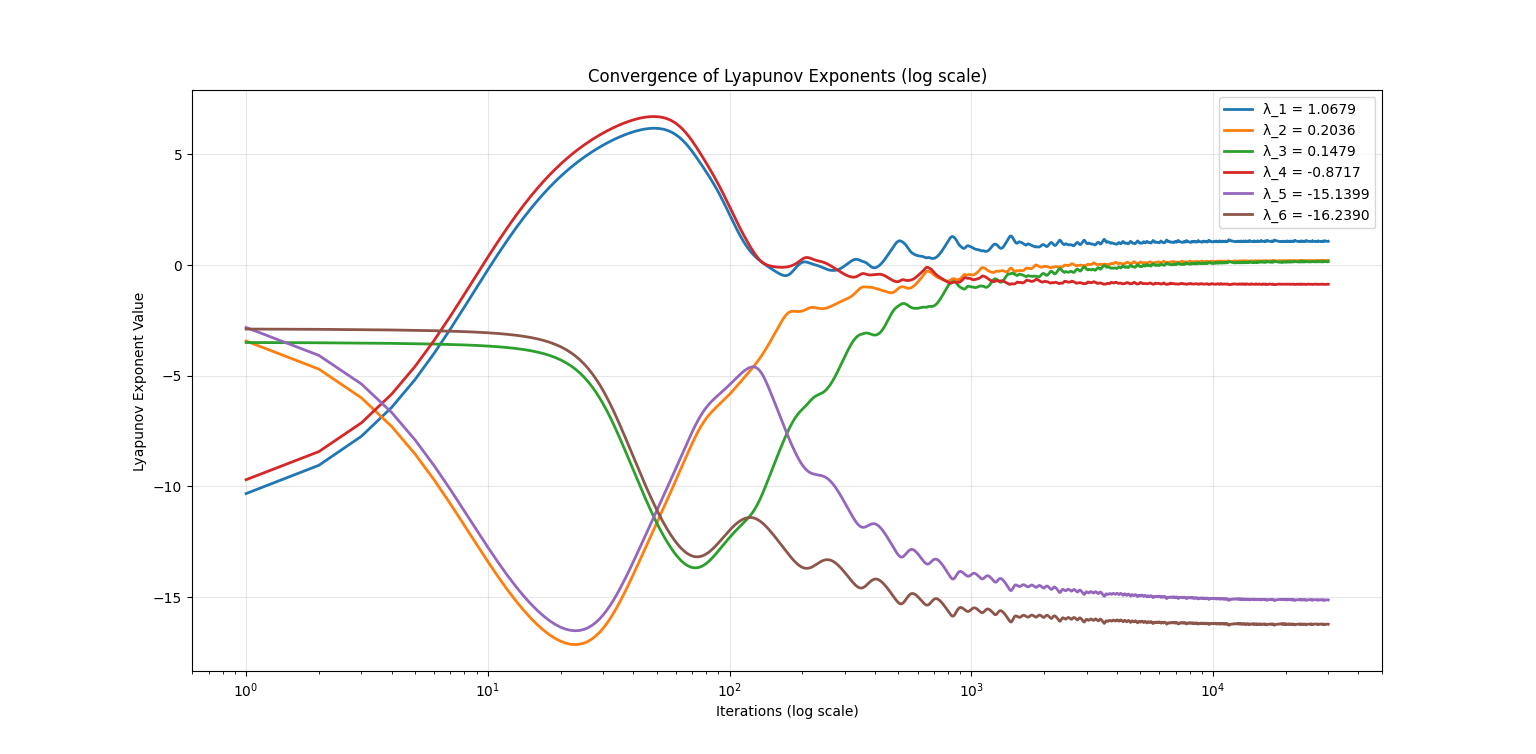}
\caption{Convergence of Lyapunov exponents over 30,000 iterations, demonstrating stabilization after approximately 12,450 iterations.}
\label{fig:lyapunov_convergence}
\end{figure}

The presence of three positive LEs unequivocally confirms hyperchaotic behavior. The LE sum is $-30.8312$, consistent with dissipative nature ($\nabla \cdot \mathbf{F} < 0$).
\\
The convergence of Lyapunov exponents over 30,000 iterations is depicted in Figure \eqref{fig:lyapunov_convergence}, demonstrating stabilization after approximately 12,450 iterations.

\subsubsection{Kaplan-Yorke Dimension}
The Kaplan-Yorke dimension $D_{KY}$ estimates the fractal dimension of the attractor from the Lyapunov spectrum:
\begin{equation*}
D_{KY} = k + \frac{\sum_{i=1}^{k} \lambda_i}{|\lambda_{k+1}|}
\end{equation*}
where $k$ is the largest integer satisfying $\sum_{i=1}^{k} \lambda_i \geq 0$.
\\
For our spectrum:
\begin{align*}
\sum_{i=1}^{1} \lambda_i &= 1.0679 > 0, \\
\sum_{i=1}^{2} \lambda_i &= 1.0679 + 0.2036 = 1.2715 > 0, \\
\sum_{i=1}^{3} \lambda_i &= 1.2715 + 0.1479 = 1.4194 > 0, \\
\sum_{i=1}^{4} \lambda_i &= 1.4194 - 0.8717 = 0.5477 > 0, \\
\sum_{i=1}^{5} \lambda_i &= 0.5477 - 15.1399 = -14.5922 < 0.
\end{align*}
Thus, $k = 4$ and:
\begin{equation*}
D_{KY} = 4 + \frac{0.5477}{15.1399} = 4.0362
\end{equation*}
This value substantially exceeds $D_{KY} \approx 2.06$ for a single Lorenz system, confirming enhanced complexity of our coupled hyperchaotic attractor, as quantified in Table \ref{table:complexity_comparison}.

\begin{table}[!ht]
\centering
\caption{Comparison of dynamical complexity metrics between classical Lorenz system and proposed coupled system}
\label{table:complexity_comparison}
\begin{tabular}{@{}lcccc@{}}
\toprule
\textbf{System} & \textbf{Positive LEs} & \textbf{$D_{KY}$} & \textbf{Dimension} & \textbf{Behavior} \\
\midrule
Classical Lorenz (2025) \cite{Bakri2025} & 1 & 2.062 & 3 & Chaotic \\
Proposed System & 3 & 4.036 & 6 & Hyperchaotic \\
\bottomrule
\end{tabular}
\end{table}

\section{Secure Communication Framework Based on the Proposed Coupled Hyperchaotic System}\label{sec5}

To illustrate the practical applicability of the proposed coupled Lorenz system, a secure communication framework based on chaos synchronization is considered. The objective of this section is not to introduce a new cryptographic architecture, but rather to demonstrate how the hyperchaotic dynamics and synchronization properties established in the previous sections can be exploited in secure transmission applications.

\begin{figure}[H]
\centering
\includegraphics[width=15cm,height=6cm]{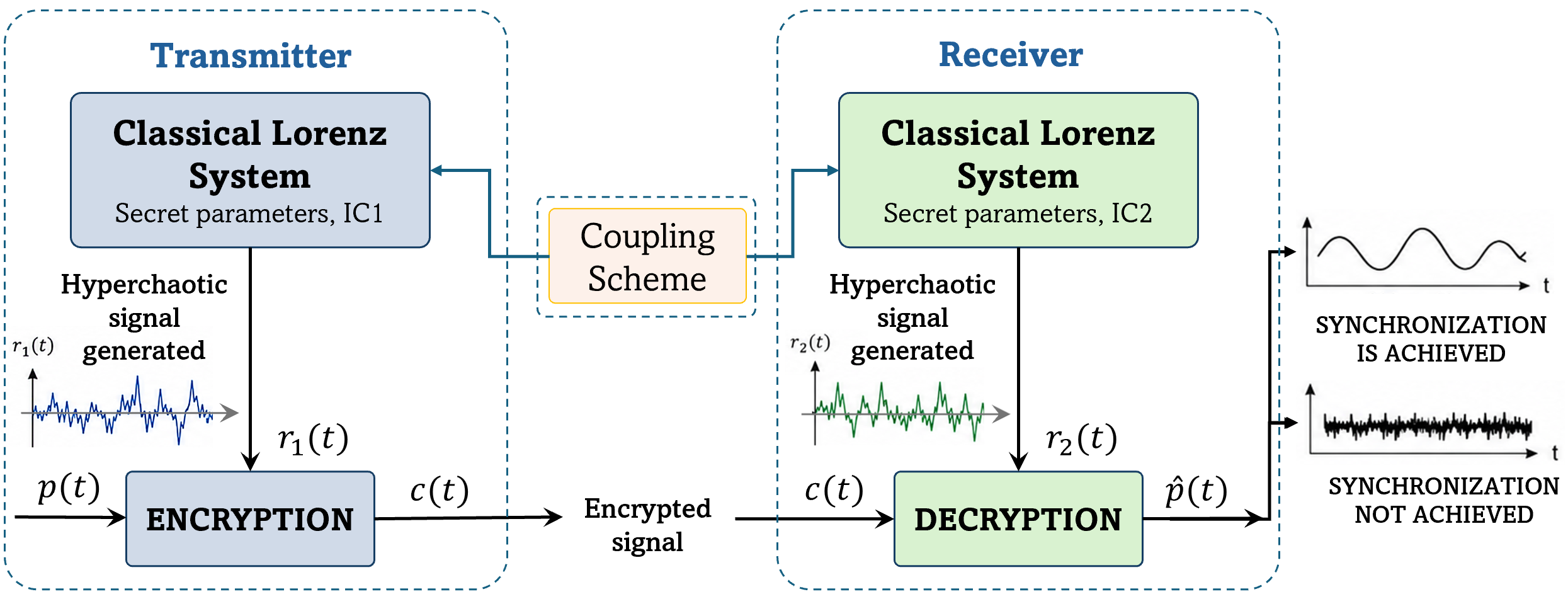}
\caption{High-level secure communication framework based on the proposed synchronized hyperchaotic coupled Lorenz system.}
\label{fig-secure-comm}
\end{figure}

As illustrated in Figure \ref{fig-secure-comm}, the transmitter and receiver each contain a chaotic generator governed by the proposed coupled Lorenz system. Initially, the generated chaotic signals are not synchronized. Through the bidirectional coupling mechanism, the receiver progressively synchronizes with the transmitter, and both chaotic generators eventually evolve along the same trajectory.
\\
At the transmitter side, the plaintext signal $p(t)$ is combined with the chaotic sequence generated by the hyperchaotic system to produce the encrypted signal $c(t)$. This encrypted signal is then transmitted through the communication channel. At the receiver side, the received signal $c(t)$ is processed using the locally generated chaotic sequence. Successful recovery of the original information is only possible after synchronization between the transmitter and receiver has been achieved.
\\
Before synchronization is established, the chaotic sequences generated by the transmitter and receiver differ significantly due to the sensitive dependence on initial conditions $(IC_1, IC_2)$. Consequently, the receiver cannot reconstruct the original message, and the recovered output remains unintelligible. Once synchronization is attained, both systems generate identical chaotic sequences, enabling accurate decryption and faithful recovery of the transmitted plaintext.
\\
The proposed coupled hyperchaotic system is particularly attractive for secure communication applications. Owing to the presence of multiple positive Lyapunov exponents, the generated signals exhibit enhanced dynamical complexity, stronger sensitivity to secret parameters and initial conditions, and larger effective key spaces compared with conventional chaotic systems. These properties contribute to improving the security level of chaos-based communication schemes while maintaining reliable synchronization performance between the transmitter and receiver.

\begin{figure}[H]
\centering
\includegraphics[width=16cm,height=12.87cm]{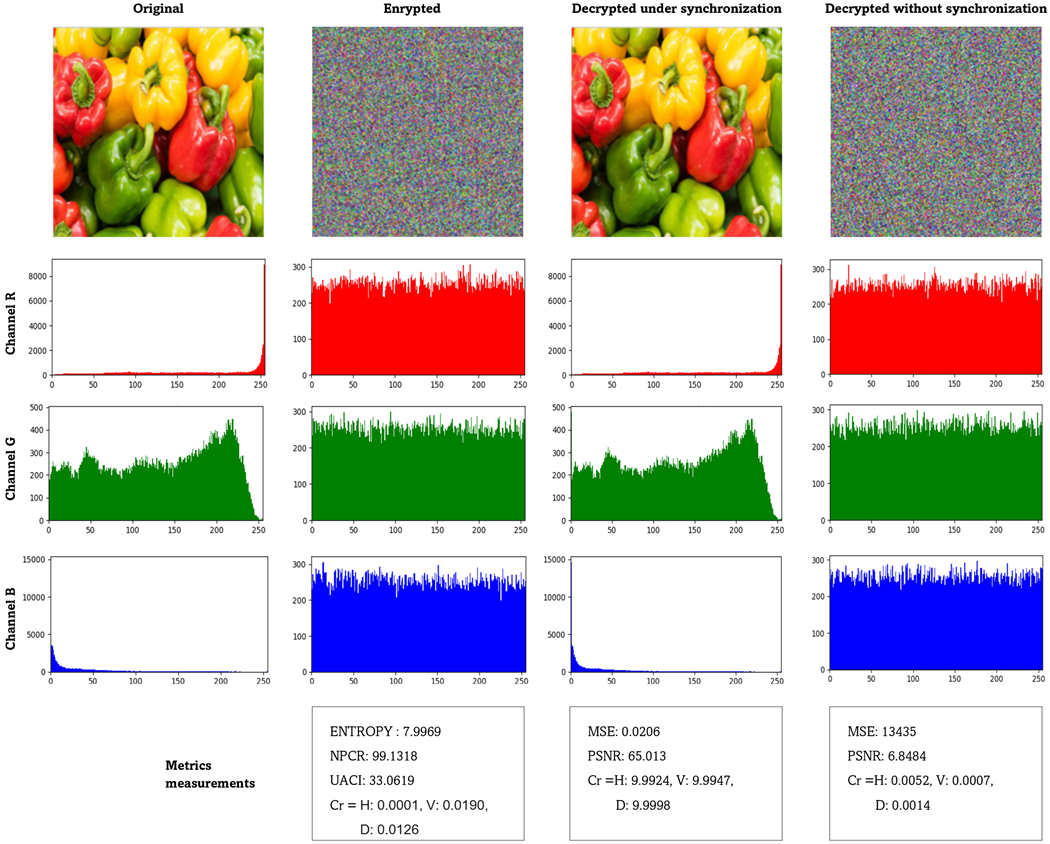}
\caption{Security and reconstruction metrics (Entropy, NPCR, UACI, MSE, PSNR) for R, G, B channels under synchronized and unsynchronized conditions.}
\label{Results}
\end{figure}

To further validate the applicability of the proposed synchronized hyperchaotic system in secure communications, an image encryption and decryption experiment was performed. The first three columns of Figure \ref{Results} show the original image, the encrypted image, and the decrypted image obtained under synchronization, together with their corresponding histograms and performance metrics. The experiments were carried out using the synchronization parameters $(\alpha_1 = 12)$ and $(\alpha_2 = 6)$.

The obtained results demonstrate that the encrypted image is visually indistinguishable from the original one and exhibits a nearly uniform histogram distribution, indicating an effective masking of the statistical characteristics of the plaintext image. Moreover, the entropy of the encrypted image approaches the theoretical ideal value, confirming the high randomness of the generated ciphertext and the strong confusion capability provided by the hyperchaotic sequences.
At the receiver side, the decrypted image is successfully reconstructed with high fidelity, demonstrating the effectiveness of the proposed synchronization mechanism. Small discrepancies may be observed during the initial transient phase due to the synchronization process, where the chaotic sequences generated at the transmitter and receiver are not yet perfectly matched. However, owing to the rapid synchronization achieved by the proposed coupling scheme, these differences rapidly vanish and have a negligible impact on the recovered image quality.

As illustrated in the fourth column of Figure \ref{Results}, the decryption process performed under unsynchronized conditions, using $(\alpha_1 = 0.5)$ and $(\alpha_2 = 0.2)$, fails to recover the original image correctly. Since the receiver chaotic generator is not synchronized with the transmitter, the generated chaotic sequences differ significantly from those employed during encryption. As a result, the decrypted image remains highly distorted and visually unintelligible. The corresponding histograms, together with the quantitative performance metrics, are also reported in the fourth column of Figure \ref{Results}, confirming the inability of the receiver to reconstruct the original image in the absence of synchronization.
\section{Comparative Analysis of Synchronization Performance}\label{sec6}

It should be emphasized that the condition $\alpha_1+\alpha_2>15.1$ constitutes a sufficient synchronization condition derived from the Lyapunov framework and Sylvester's criterion. Therefore, it guarantees the asymptotic convergence of the synchronization error but does not represent a necessary condition.

Consequently, synchronization may still occur for smaller values of the coupling parameter. This behavior is a well-known consequence of the conservatism inherent to Lyapunov-based sufficient criteria. To illustrate this point and assess the conservatism of the theoretical bound, additional numerical investigations are presented in the following section. Several representative cases with ($\alpha_1+\alpha_2<15.1$) are examined, showing that complete synchronization can still be achieved despite the fact that the sufficient analytical conditions are not fulfilled.

\begin{figure}[H]
\centering
\includegraphics[width=15cm,height=6cm]{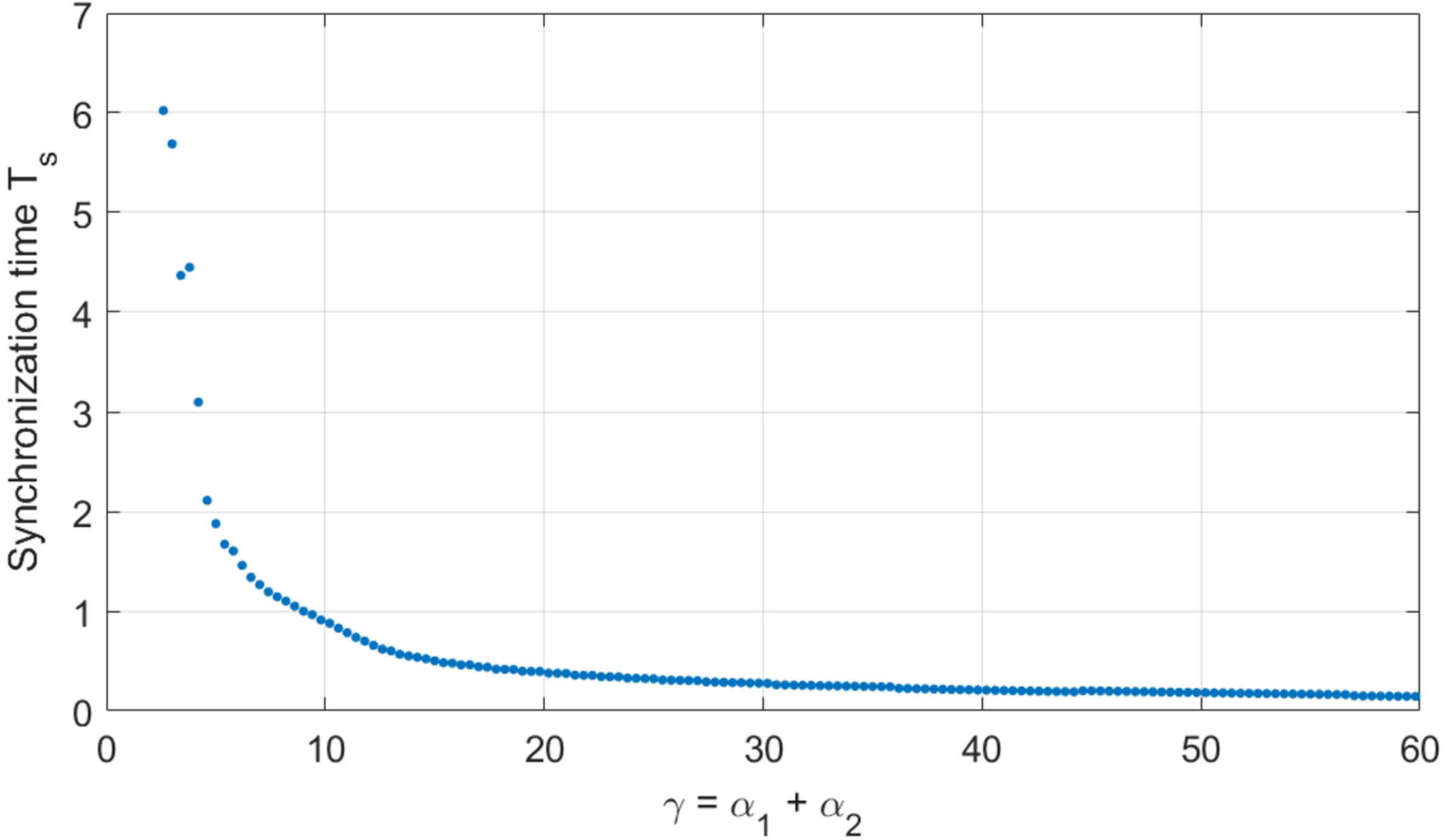}
\caption{Synchronization time versus the coupling parameter $\gamma=\alpha_1+\alpha_2$.}\label{fig:A}
\end{figure}

Compared with existing synchronization methods for coupled Lorenz systems, the proposed approach presents three main advantages: (i) a bidirectional asymmetric coupling structure, (ii) an explicit analytical synchronization condition expressed directly in terms of the coupling coefficients, and (iii) improved synchronization performance in terms of convergence speed. To evaluate this latter aspect, the synchronization time was computed as a function of the coupling parameter $\gamma=\alpha_1+\alpha_2$. The obtained results (presented by Figure \ref{fig:A}) show that for relatively weak coupling strengths ($\gamma<5$), the synchronization time remains close to those reported in the literature, namely about 4.5 s in \cite{Ogabi2025}, 5 s in \cite{Chauhan2023}, and 5.4 s in \cite{Akter2023}, and may even become larger for smaller values of $\gamma$. However, as the coupling strength increases beyond $\gamma=5$, the synchronization time decreases rapidly and continuously, approaching values close to zero for sufficiently large $\gamma$. This behavior demonstrates that the proposed coupling scheme not only guarantees synchronization under explicit analytical conditions but also significantly accelerates the convergence process compared with conventional synchronization approaches reported in the literature.

\section*{Conclusion}

This paper presented a six-dimensional hyperchaotic system obtained through the bidirectional asymmetric diffusive coupling of two classical Lorenz oscillators. A comprehensive dynamical analysis, including dissipativity, equilibrium-point stability, Lyapunov spectrum, bifurcation, and synchronization analyses, was carried out to characterize the proposed system. The results demonstrated the existence of broad hyperchaotic operating regions with up to three positive Lyapunov exponents, providing a practical framework for selecting the coupling parameters. In particular, the proposed analysis shows that choosing a coupling strength $\gamma = \alpha_1 + \alpha_2 > 15.1$ simultaneously guarantees system dissipativity, robust hyperchaotic behavior, and complete synchronization, thereby providing a simple and effective guideline for practical implementation. 

A Lyapunov-based synchronization analysis led to explicit sufficient synchronization conditions together with a less conservative synchronization criterion based on effective attractor bounds. Numerical simulations confirmed the theoretical analysis and showed that the synchronization time decreases significantly as the coupling strength increases, reaching values below $0.5$ s, compared with approximately $4.5$--$5.4$ s reported by existing synchronization approaches. The bifurcation analysis further validated the parameter regions identified from the Lyapunov spectrum and confirmed the persistence of chaotic dynamics over a wide range of coupling strengths. Finally, image encryption and decryption experiments successfully demonstrated that the proposed synchronized hyperchaotic system enables reliable secure communication while benefiting from the increased complexity, enlarged key space, and enhanced security offered by hyperchaotic dynamics.

\vspace{0.5cm}

\section*{Acknowledgments}
This research was supported by the Laboratory of Defense Sciences and Technologies (LR19DN01) of the Tunisian Ministry of National Defense.


\begin{thebibliography}{00}

\bibitem{Akter2023}
Akter, M. T., Tarammim, A., and Hussen, S. (2023).
\emph{Chaos control and synchronization of modified Lorenz system using active control and backstepping scheme}.
Waves in Random and Complex Media, 1--20.

\bibitem{Bakri2025}
Bakri, T., and Verhulst, F. (2025).
\emph{A note on the Kaplan-Yorke dimension}.
International Journal of Bifurcation and Chaos, 35(1), 2550133.

\bibitem{Biswas2026} Biswas, F., Biswas, S., Islam Mondal, R., Goldar, S., Sarif Hassan, S., and Sardar, P. (2026).
\emph{Chaos and Complexity in a Novel Three-Dimensional Nonlinear System: Analysis, Bifurcations, and Applications}.
Mathematical Methods in the Applied Sciences.

\bibitem{bragard2007}
Bragard, J., Vidal, G., Mancini, S., Mendoza, C., and Boccaletti, S. (2007).
\emph{Synchronization of coupled Lorenz systems}.
Physical Review E, 75(2), 026211.

\bibitem{9s}
Brown, R., and Kocarev, L. (2000).
\emph{A unifying definition of synchronization for dynamical systems}.
Chaos, 10(2), 344--349.

\bibitem{chen2004}
Chen, A., Lu, J., Lü, J., and Yu, S. (2006).
\emph{Generating hyperchaotic Lü attractor via state feedback control}.
Physica A: Statistical Mechanics and its Applications, 364, 103--110.

\bibitem{CU1999}
Chen, G., and Ueta, T. (1999).
\emph{Yet another chaotic attractor}.
International Journal of Bifurcation and Chaos, 9(7), 1465--1466.

\bibitem{Chauhan2023}
Chauhan, Y., and Sharma, B. (2023).
\emph{Synchronisation results for an interconnected network of nonlinear systems with diffusive nonlinear coupling using contraction}.
Journal of Applied Nonlinear Dynamics, 12(3), 591.

\bibitem{1s}
Chiu, C., Lin, W., and Peng, C. (2000).
\emph{Asymptotic synchronization in lattices of coupled nonlinear Lorenz equations}.
International Journal of Bifurcation and Chaos, 10(12), 2717--2728.

\bibitem{5s}
Chua, L. O., Itah, M., Kosarev, L., and Eckert, K. (1993).
\emph{Chaos synchronization in Chua's circuits}.
Journal of Circuits, Systems and Computers, 3(1), 93--108.

\bibitem{gao2008}
Gao, T., and Chen, Z. (2008).
\emph{A new image encryption algorithm based on hyper-chaos}.
Physics Letters A, 372(4), 394--400.

\bibitem{grassi2002}
Grassi, G., and Miller, D. A. (2002).
\emph{Theory and experimental realization of observer-based discrete-time hyperchaos synchronization}.
IEEE Transactions on Circuits and Systems I: Fundamental Theory and Applications, 49(3), 373--378.

\bibitem{grassi2008}
Grassi, G., Severance, F. L., Mashev, D. L., Bazuin, B. J., and Miller, D. A. (2008).
\emph{Four-wing attractor in symmetrically coupled Lorenz systems}.
IEEE Transactions on Circuits and Systems II: Express Briefs, 55(11), 1194--1198.

\bibitem{6s}
Hu, G., Xiao, J., and Zheng, Z. (2000).
\emph{Chaos control}.
Shanghai Scientific and Technological Education Publishing House.

\bibitem{kanso2012}
Kanso, A., and Ghebleh, M. (2012).
\emph{A novel image encryption algorithm based on a 3D chaotic map}.
Communications in Nonlinear Science and Numerical Simulation, 17(7), 2943--2959.

\bibitem{KS2015}
Khan, A., and Singh, P. (2015).
\emph{Chaos synchronization in Lorenz system}.
Applied Mathematics, 6(11), 1864--1872.

\bibitem{li2005}
Li, S., Chen, G., and Mou, X. (2005).
\emph{On the dynamical degradation of digital piecewise linear chaotic maps}.
International Journal of Bifurcation and Chaos, 15(10), 3119--3151.

\bibitem{lorenz1963}
Lorenz, E. N. (1963).
\emph{Deterministic nonperiodic flow}.
Journal of the Atmospheric Sciences, 20(2), 130--141.

\bibitem{13s}
Lü, J., Lu, J., and Chen, S. (2002).
\emph{Chaotic time series analysis and its application}.
Wuhan University Press.

\bibitem{Ogabi2025}
Ogabi, C., Shehu, T., Idowu, B., Fakunle, O., Onori, E., Mustapha, R., Bamgbose, M., and Muhammad, S. (2025).
\emph{Stability, Tracking Control, and Synchronisation of Hyperchaotic 5D Lorenz System using Active Backstepping Designs}.
Journal of Research and Review in Science, 11(2), 28--48.

\bibitem{Garcia2018}
Ontañón-García, L. J., Martínez, M. G., Cantón, I. C., Montalvo, C. S., Torres, M. R., and Ponce, R. L. (2018).
\emph{Bifurcation from chaos to periodic states in bidirectional interconnected Lorenz systems by the variation of the coupling strengths}.
IFAC-PapersOnLine, 51(33), 86--90.

\bibitem{pecora1990}
Pecora, L. M., and Carroll, T. L. (1990).
\emph{Synchronization in chaotic systems}.
Physical Review Letters, 64(8), 821--824.

\bibitem{wang2008}
Wang, X., and Wang, M. (2008).
\emph{A hyperchaos generated from Lorenz system}.
Physica A: Statistical Mechanics and its Applications, 387(14), 3751--3758.

\bibitem{2s}
Wang, C., and Ge, S. S. (2001).
\emph{Adaptive synchronization of uncertain chaotic systems via backstepping design}.
Chaos, Solitons and Fractals, 12(7), 1199--1206.

\bibitem{Y2007}
Yu, Y. (2007).
\emph{The synchronization for time-delay of linearly bidirectional coupled chaotic system}.
Chaos, Solitons \& Fractals, 33(4), 1197--1203.

\bibitem{YZ2004}
Yu, Y., and Zhang, S. (2004).
\emph{The synchronization of linearly bidirectional coupled chaotic systems}.
Chaos, Solitons \& Fractals, 22(1), 189--197.

\bibitem{YZZZ2022}
Yu, X., Zuo, Z., Zhu, S., and Zhang, X. (2022).
\emph{A General Criterion for Unidirectionally Coupled Generalized Chaotic Synchronization with a Desired Manifold}.
Advances in Mathematical Physics, 2022(1), 5088434.

\bibitem{Zhou2026} Zhou, Y., Guo, Y., Zhang, W., and Zhao, P. (2026).
\emph{Synchronization, circuit implementation, and a three-layer image encryption scheme for a novel 5D hyperchaotic system}.
Nonlinear Dynamics, 114(6), 416.
\end{thebibliography}
\end{document}